\documentclass[a4paper,10pt]{amsart}

\usepackage{amsmath,amssymb,amsthm}

\newcommand{\CM}{\overline{\operatorname{\mathcal{M}}}}

\newcommand{\IM}{\operatorname{\mathcal{M}}}
\newcommand{\BB}{\operatorname{\mathcal{B}}}
\newcommand{\EE}{\operatorname{\mathcal{E}}}
\newcommand{\LL}{\operatorname{\mathcal{L}}}

\renewcommand{\AA}{\operatorname{\mathfrak{A}}}
\newcommand{\PP}{\operatorname{\mathfrak{P}}}
\newcommand{\WW}{\operatorname{\mathfrak{W}}}
\newcommand{\CC}{\operatorname{\mathfrak{C}}}
\newcommand{\DD}{\operatorname{\mathfrak{D}}}

\newcommand{\Si}{\dot{S}}

\newcommand{\CP}{\mathbb{C}\mathbb{P}^1}

\newcommand{\CR}{\bar{\partial}}

\newcommand{\SFT}{\operatorname{SFT}}
\newcommand{\str}{\operatorname{string}}
\newcommand{\BV}{\operatorname{BV}}
\newcommand{\ev}{\operatorname{ev}}

\newcommand{\cst}{\operatorname{const}}

\newcommand{\CZ}{\operatorname{CZ}}
\newcommand{\Morse}{\operatorname{Morse}}

\newcommand{\rank}{\operatorname{rank}}

\newcommand{\coker}{\operatorname{coker}}
\newcommand{\del}{\partial}

\newcommand{\RS}{\IR \times S^1}
\newcommand{\Ju}{\underline{J}}

\renewcommand{\qed}{\square}

\newcommand{\im}{\operatorname{im}}
\newcommand{\IC}{\operatorname{\mathbb{C}}}
\newcommand{\IZ}{\operatorname{\mathbb{Z}}}

\newcommand{\IR}{\operatorname{\mathbb{R}}}
\newcommand{\IN}{\operatorname{\mathbb{N}}}
\newcommand{\ID}{\operatorname{\mathbf{D}}}
\newcommand{\IH}{\operatorname{\mathbf{H}}}
\newcommand{\IG}{\operatorname{\mathbf{G}}}
\newcommand{\IL}{\operatorname{\mathbf{L}}}
\newcommand{\Ih}{\operatorname{\mathbf{h}}}
\newcommand{\Ig}{\operatorname{\mathbf{g}}}
\newcommand{\Il}{\operatorname{\mathbf{l}}}
\newcommand{\IF}{\operatorname{\mathbf{F}}}

\newcommand{\Coker}{\operatorname{Coker}}

\newcommand{\kom}[1]{}

\date{March 2010}
\author{Oliver Fabert}
\title{Gravitational descendants in\\ symplectic field theory}
\thanks{Research supported by the German Research Foundation (DFG)}

\begin{document}

\abstract 
It was pointed out by Y. Eliashberg in his ICM 2006 plenary talk that the rich algebraic formalism of symplectic field theory leads to a natural appearance of quantum and classical integrable systems, at least in the case when the contact manifold is the prequantization space of a symplectic manifold. In this paper we generalize the definition of gravitational descendants in SFT from circle bundles in the Morse-Bott case to general contact manifolds. After we have shown using the ideas in [OP] that for the basic examples of holomorphic curves in SFT, that is, branched covers of cylinders over closed Reeb orbits, the gravitational descendants have a geometric interpretation in terms of branching conditions, we follow the ideas in [CL] to compute the corresponding sequence of Poisson-commuting functions when the contact manifold is the unit cotangent bundle of a Riemannian manifold.}

\maketitle

\tableofcontents

\markboth{OLIVER FABERT}{DESCENDANTS IN SFT}

\section*{Summary}

Symplectic field theory (SFT), introduced by H. Hofer, A. Givental and Y. Eliashberg in 2000 ([EGH]), is a very large project and can be viewed as a topological quantum field theory approach to Gromov-Witten theory. Besides providing a unified view on established pseudoholomorphic curve theories like symplectic Floer homology, contact homology and Gromov-Witten theory, it leads to numerous new applications and opens new routes yet to be explored. \\
 
While symplectic field theory leads to algebraic invariants with very rich algebraic structures, which are currently studied by a large group of researchers, for all the geometric applications found so far it was sufficient to work with simpler invariants like cylindrical contact homology. Although cylindrical contact homology is not always defined, it is much easier to compute, not only since it involves just moduli spaces of holomorphic cylinders but also due to the simpler algebraic formalism. While the rich algebraic formalism of the higher invariants of symplectic field theory seems to be too complicated for concrete geometric applications, it was pointed out by Eliashberg in his ICM 2006 plenary talk ([E]) that the integrable systems of rational Gromov-Witten theory very naturally appear in rational symplectic field theory by using the link between the rational symplectic field theory of circle bundles in the Morse-Bott version and the rational Gromov-Witten potential of the underlying symplectic manifold. Indeed, after introducing gravitational descendants as in Gromov-Witten theory, it is precisely the rich algebraic formalism of SFT with its Weyl and Poisson structures that provides a natural link between symplectic field theory and (quantum) integrable systems. In particular, in the case where the contact manifold is a circle bundle over a closed symplectic manifold, the rich algebraic formalism of symplectic field theory seems to provide the right framework to understand the deep relation between Gromov-Witten theory and integrable systems, at least in the genus zero case. \\

While in the Morse-Bott case in [E] it follows from the corresponding statements for the Gromov-Witten descendant potential that the sequences of commuting operators and Poisson-commuting functions are independent of auxiliary choices like almost complex structure and abstract perturbations, for the case of general contact manifolds it is well-known that the SFT Hamiltonian however in general explicitly depend on choices like contact form, cylindrical almost complex structure and coherent abstract perturbations and hence is not an invariant for the contact manifold itself. But before we can come down to the question of invariance, we first need to give a rigorous definition of gravitational descendants in the context of symplectic field theory. \\

While in Gromov-Witten theory the gravitational descendants were defined by integrating powers of the first Chern class of the tautological line bundle over the moduli space, which by Poincare duality corresponds to counting common zeroes of sections in this bundle, in symplectic field theory, more generally every holomorphic curves theory where curves with punctures and/or boundary are considered, we are faced with the problem that the moduli spaces generically have codimension-one boundary, so that the count of zeroes of sections in general depends on the chosen sections in the boundary. It follows that the integration of the first Chern class of the tautological line bundle over a single moduli space has to be replaced by a construction involving all moduli space at once. Note that this is similar to the choice of coherent abstract perturbations for the moduli spaces in symplectic field theory in order to achieve transversality for the Cauchy-Riemann operator. Keeping the interpretation of descendants as common zero sets of sections in powers of the tautological line bundles (which will turn out to be particularly useful when one studies the topological meaning of descendants by localizing on special divisors, see [FR]), we define in this paper the notion of {\it coherent collections of sections} in the tautological line bundles over all moduli spaces, which just formalizes how the sections chosen for the lower-dimensional moduli spaces should affect the section chosen for a moduli spaces on its boundary. To be more precise, since the sections should be invariant under obvious symmetries like reordering of the punctures and the marked points, we actually need to work with multi-sections in order to meet both the symmetry and the transversality assumption. We will then define {\it descendants of moduli spaces} $\CM^j\subset\CM$, which we obtain inductively as zero sets of these coherent collections of sections $(s_j)$ in the tautological line bundles over the descendant moduli spaces $\CM^{j-1}\subset\CM$, and define {\it descendant Hamiltonians} $\IH^1_{i,j}$ by integrating chosen closed differential forms $\theta_i$ over $\CM^j$. For these we prove the following theorem. \\ 
\\
{\bf Theorem:} {\it Counting holomorphic curves with one marked point after integrating differential forms and introducing gravitational descendants defines a sequence of distinguished elements} 
\begin{equation*} \IH^1_{i,j}\in H_*(\hbar^{-1}\WW^0,D^0) \end{equation*} 
{\it in the full SFT homology algebra with differential $D^0=[\IH^0,\cdot]: \hbar^{-1}\WW^0\to\hbar^{-1}\WW^0$, which commute with respect to the commutator bracket on $H_*(\hbar^{-1}\WW^0,D^0)$,} 
\begin{equation*} [\IH^1_{i,j},\IH^1_{k,\ell}] = 0,\; (i,j),(k,\ell)\in\{1,...,N\}\times\IN. \end{equation*}
\\
In contrast to the Morse-Bott case considered in [E] it follows that, when the differential in symplectic field theory counting holomorphic curves without additional marked points is no longer zero, the sequences of generating functions no longer commute with respect to the bracket, but only commute {\it after passing to homology.} On the other hand, in the same way as the rational symplectic field theory of a contact manifold is defined by counting only curves with genus zero, we immediately obtain a rational version of the above statement by expanding $\IH^0$ and the $\IH^1_{i,j}$ in powers of the formal variable $\hbar$ for the genus. \\
\\
{\bf Corollary:} {\it Counting rational holomorphic curves with one marked point after integrating differential forms and introducing gravitational descendants defines a sequence of distinguished elements} 
\begin{equation*} \Ih^1_{i,j}\in H_*(\PP^0,d^0), \end{equation*} 
{\it in the rational SFT homology algebra with differential $d^0=\{\Ih^0,\cdot\}: \PP^0\to\PP^0$, which commute with respect to the Poisson bracket on $H_*(\PP^0,d^0)$,} 
\begin{equation*} \{\Ih^1_{i,j},\Ih^1_{k,\ell}\} = 0,\; (i,j),(k,\ell)\in\{1,...,N\}\times\IN. \end{equation*}
\\
As we already outlined above, in contrast to the circle bundle  case we have to expect that the sequence of descendant Hamiltonians depends on the auxiliary choices like contact form, cylindrical almost complex structure and coherent abstract polyfold perturbations. Here we prove the following natural invariance statements. \\
\\
{\bf Theorem:} {\it For different choices of contact form $\lambda^{\pm}$, cylindrical almost complex structure $\Ju^{\pm}$ , abstract polyfold perturbations and sequences of coherent collections of sections $(s^{\pm}_j)$ the resulting systems of commuting operators $\IH^{1,-}_{i,j}$ on $H_*(\hbar^{-1}\WW^{0,-},D^{0,-})$ and $\IH^{1,+}_{i,j}$ on $H_*(\hbar^{-1}\WW^{0,+},D^{0,+})$ are isomorphic, i.e., there exists an isomorphism of the Weyl algebras $H_*(\hbar^{-1}\WW^{0,-},D^{0,-})$ and $H_*(\hbar^{-1}\WW^{0,+},D^{0,+})$ which maps $\IH^{1,-}_{i,j}\in H_*(\hbar^{-1}\WW^{0,-},D^{0,-})$ to $\IH^{1,+}_{i,j}\in H_*(\hbar^{-1}\WW^{0,+},D^{0,+})$.}\\ 
\\
Note that this theorem is an extension of the theorem in [EGH] stating that for different choices of auxiliary data the Weyl algebras $H_*(\hbar^{-1}\WW^{0,-},D^{0,-})$ and $H_*(\hbar^{-1}\WW^{0,+},D^{0,+})$ are isomorphic. As above we clearly also get a rational version of the invariance statement: \\
\\
{\bf Corollary:} {\it For different choices of contact form $\lambda^{\pm}$, cylindrical almost complex structure $\Ju^{\pm}$, abstract polyfold perturbations and sequences of coherent collections of sections $(s^{\pm}_j)$ the resulting system of Poisson-commuting functions $\Ih^{1,-}_{i,j}$ on $H_*(\PP^{0,-},d^{0,-})$ and $\Ih^{1,+}_{i,j}$ on $H_*(\PP^{0,+},d^{0,+})$ are isomorphic, i.e., there exists an isomorphism of the Poisson algebras $H_*(\PP^{0,-},d^{0,-})$ and $H_*(\PP^{0,+},d^{0,+})$ which maps $\Ih^{1,-}_{i,j}\in H_*(\PP^{0,-},d^{0,-})$ to $\Ih^{1,+}_{i,j}\in H_*(\PP^{0,+},d^{0,+})$.} \\

As concrete example beyond the case of circle bundles discussed in [E] we consider the symplectic field theory of a closed geodesic. For this recall that in [F2] the author introduces the symplectic field theory of a closed Reeb orbit $\gamma$, which is defined by counting only those holomorphic curves which are branched covers of the orbit cylinder $\IR\times\gamma$ in $\IR\times V$. In [F2] we prove that these orbit curves do not contribute to the algebraic invariants of symplectic field theory as long as they do not carry additional marked points. Our proof explicitly uses that the subset of orbit curves over a fixed orbit is closed under taking boundaries and gluing, which follows from the fact that they are also trivial in the sense that they have trivial contact area and that this contact area is preserved under taking boundaries and gluing. It follows that every algebraic invariant of symplectic field theory has a natural analog defined by counting only orbit curves. In particular, in the same way as we define sequences of descendant Hamiltonians $\IH^1_{i,j}$ and $\Ih^1_{i,j}$ by counting general curves in the symplectization of a contact manifold, we can define sequences of descendant Hamiltonians $\IH^1_{\gamma,i,j}$ and $\Ih^1_{\gamma,i,j}$ by just counting branched covers of the orbit cylinder over $\gamma$ with signs (and weights), where the preservation of the contact area under splitting and gluing of curves proves that for every theorem from above we have a version for $\gamma$. We further prove that for branched covers of orbit cylinders over any closed Reeb orbit the gravitational descendants indeed have a geometric interpretation in terms of branching conditions, which generalizes the work of [OP] used in [E] for the circle. \\

Since all the considered holomorphic curves factor through the embedding of the closed Reeb orbit into the contact manifold, it follows that it only makes sense to consider differential forms of degree zero or one.  While it follows from the result $\Ih^0_{\gamma}=0$ in [F2] that the sequences $\Ih^1_{\gamma,i,j}$ indeed commute with respect to the Poisson bracket (before passing to homology), the same proof as in [F2] shows that every descendant Hamiltonian in the sequence vanishes if the differential form is of degree zero. For differential forms of degree one the strategy of the proof however no longer applies and it is indeed shown in [E] that for $\gamma=V=S^1$ and $\theta=dt$ we get nontrivial contributions from branched covers. In this paper we want to determine the corresponding Poisson-commuting sequence in the special case where the contact manifold is the unit cotangent bundle $S^*Q$ of a ($m$-dimensional) Riemannian manifold $Q$, so that every closed Reeb orbit $\gamma$ on $V=S^*Q$ corresponds to a closed geodesic $\bar{\gamma}$ on $Q$. For this we denote by $\WW^0_{\gamma}$ be the graded Weyl subalgebra of the Weyl algebra $\WW^0$, which is generated only by those $p$- and $q$-variables $p_n=p_{\gamma^n}$, $q_n=q_{\gamma^n}$ corresponding to Reeb orbits which are multiple covers of the fixed orbit $\gamma$ {\it and which are good in the sense of [BM]}. In the same way we further introduce the Poisson subalgebra $\PP^0_{\gamma}$ of $\PP^0$. Setting $q_{-n}=p_n$ we prove the following \\
\\
{\bf Theorem:} {\it Assume that the contact manifold is the unit cotangent bundle $V=S^*Q$ of a Riemannian manifold $Q$, so that the closed Reeb orbit $\gamma$ corresponds to a closed geodesic $\bar{\gamma}$ on $Q$, and that the string of differential forms just consists of a single one-form which integrates to one around the orbit. Then the resulting system of Poisson-commuting functions $\Ih^1_{\gamma,j}$, $j\in\IN$ on $\PP^0_{\gamma}$ is isomorphic to the system of Poisson-commuting functions $\Ig^1_{\bar{\gamma},j}$, $j\in\IN$ on $\PP^0_{\bar{\gamma}}=\PP^0_{\gamma}$, where for every $j\in\IN$ the descendant Hamiltonian $\Ig^1_{\bar{\gamma},j}$ is given by} 
\begin{equation*} 
 \Ig^1_{\bar{\gamma},j} \;=\; \sum \epsilon(\vec{n})\;\frac{q_{n_1}\cdot ... \cdot q_{n_{j+2}}}{(j+2)!} 
\end{equation*}
{\it where the sum runs over all ordered monomials $q_{n_1}\cdot ... \cdot q_{n_{j+2}}$ with $n_1+...+n_{j+2} = 0$ \textbf{and which are of degree $2(m+j-3)$}. Further $\epsilon(\vec{n})\in\{-1,0,+1\}$ is fixed by a choice of coherent orientations in symplectic field theory and is zero if and only if one of the orbits $\gamma^{n_1},...,\gamma^{n_{j+2}}$ is bad.} \\

Note that in the case of the circle $\bar{\gamma}=Q=S^1$ the degree condition is automatically fulfilled and we just get back the sequence of descendant Hamiltonians for the circle in [E], which agrees with the sequence of Poisson-commuting integrals of the dispersionless KdV integrable hierarchy. Forgetting about the appearing sign issues, it follows that the sequence $\Ig^1_{\bar{\gamma},j}$ is obtained from the sequence for the circle by removing all summands with the wrong, that is, not maximal degree, so that the system is completely determined by the KdV hierarchy and the Morse indices of the closed geodesic and its iterates. \\

\noindent{\bf Remark:} Note that the signs in the formula for $\Ig^1_{\bar{\gamma},j}$ are determined by the linearized Reeb flow around $\gamma$ and a choice of 
orientations for all multiples of $\gamma$. For this recall from [BM] that in order to orient moduli spaces in symplectic field theory one additionally needs to choose orientations for all occuring Reeb orbits, while the resulting invariants are independent of these auxiliary choices. While the precise formula for the functions $\Ig^1_{\bar{\gamma},j}$ depends on these choices, the resulting systems of Poisson-commuting functions for different choices are indeed isomorphic, since changing the orientations for some orbits $\gamma^k$ leads to an automorphism of the underlying Poisson algebra. Apart from the fact that the commutativity condition $\{\Ig^1_{\bar{\gamma},j},\Ig^1_{\bar{\gamma},k}\}=0$ clearly leads to relations between the different $\epsilon(\vec{n})$, observe that a choice of orientation for $\gamma$ does {\it not} lead to a canonical choice of orientations for its multiples $\gamma^k$. While we expect that it is in general very hard to write down a set of signs $\epsilon(\vec{n})$ explicitly, for all the geometric applications we have in mind and the educational purposes as a test model beyond the Gromov-Witten case we are rather interested in proving vanishing results as the one below than giving precise formulas. \\

While in the case of the circle we obtain a complete set of integrals, it is important that our theorem allows us to prove the following vanishing result, which simply follows from the fact that for hyperbolic Reeb orbits the Conley-Zehnder index is multiplicative.  \\
\\
{\bf Corollary:} {\it Assume that the closed geodesic $\bar{\gamma}$ represents a hyperbolic Reeb orbit in the unit cotangent bundle of a surface $Q$. Then $\Ig^1_{\bar{\gamma},j}=0$ and hence $\Ih^1_{\gamma,j}=0$ for all $j>0$.} \\

Note that this result is actually true for $\dim Q > 1$. While for $\dim Q > 2$ the result directly follows from index reasons, in the case when $\dim Q=2$ a simple computation shows that all moduli spaces with $2j+1$ punctures possibly contribute to the descendant Hamiltonian $\Ih^1_{\gamma,j}$. Since in this case the Fredholm index is $2j-1$ and hence for $j>0$ strictly smaller than the dimension of the underlying nonregular moduli space of branched covers, which is $4j-2$, transversality cannot be satisfied but the cokernels of the linearized operators fit together to give an obstruction bundle of rank $2j-1$. \\

Apart from using the geometric interpretation of gravitational descendants for branched covers of orbit cylinders over a closed Reeb orbit in terms of branching conditions mentioned above, the second main ingredient for the proof is the idea in [CL] to compute the symplectic field theory of $V=S^*Q$ from the string topology of the underlying Riemannian manifold $Q$ by studying holomorphic curves in the cotangent bundle $T^*Q$. More precisely, we compute the symplectic field theory of a closed Reeb orbit $\gamma$ in $S^*Q$ including differential forms and gravitational descendants by studying branched covers of the trivial half-cylinder connecting the closed Reeb orbit in the unit cotangent bundle with the underlying closed geodesic in the cotangent bundle $T^*Q$ with special branching data, where the latter uses the geometric interpretation of gravitational descendants. In order to give a complete proof we also prove the neccessary transversality theorems using finite-dimensional obstruction bundles over the underlying nonregular moduli spaces. While on the SFT side one has very complicated obstruction bundles over nonregular moduli spaces of arbitary large dimension, on the string side all relevant nonregular moduli spaces already turn out to be discrete, so that the obstruction bundles disappear if the Fredholm index is right. It follows that the system of Poisson-commuting function for a closed geodesic is completely determined by the KdV hierarchy and the Morse indices of the closed geodesic and its iterates. \\

This paper is organized as follows. \\

Section one is concerned with the definition and the basic results about gravitational descendants in symplectic field theory. After we recalled the basic definitions of symplectic field theory in subsection 1.1, we define gravitational descendants in subsection 1.2 using the coherent collections of sections and prove that the resulting sequences of descendant Hamiltonians commute after passing to homology. In subsection 1.3 we prove the desired invariance statement and discuss the important case of circle bundles in the Morse-Bott setup outlined in [E] in 1.4. \\

After we treated the general case in section one, section two is concerned with a concrete example beyond the case of circle bundles, the symplectic field theory of a closed geodesic, which naturally generalizes the case of the circle in [E]. After we have recalled the definition of symplectic field theory for a closed Reeb orbit including the results from [F2] in subsection 2.1, we show in subsection 2.2 that for branched covers of orbit cylinders the gravitational descendants have a geometric interpretation in terms of branching conditions. After outlining that there exists a version of the isomorphism in [CL] involving the symplectic field theory of a closed Reeb orbit in the unit cotangent bundle, we study the moduli space of branched covers of the corresponding trivial half-cylinder in the cotangent bundle in subsection 2.3. Since we meet the same transversality problems as in [F2], we study the neccessary obstruction bundle setup including Banach manifolds and Banach space bundles in subsection 2.3. In subsection 2.4 we finally prove the above theorem by studying branched covers of the trivial half-cylinder with special branching behavior. \\  
\\
{\bf Acknowledgements:} This research was supported by the German Research Foundation (DFG). 
The author thanks K. Cieliebak, Y. Eliashberg, K. Fukaya, M. Hutchings and P. Rossi for useful discussions.
 
\section{Symplectic field theory with gravitational descendants}

\subsection{Symplectic field theory}
Symplectic field theory (SFT) is a very large project, initiated by Eliashberg,
Givental and Hofer in their paper [EGH], designed to describe in a unified way 
the theory of pseudoholomorphic curves in symplectic and contact topology. 
Besides providing a unified view on well-known theories like symplectic Floer 
homology and Gromov-Witten theory, it shows how to assign algebraic invariants 
to closed contact manifolds $(V,\xi=\{\lambda=0\})$: \\
   
Recall that a contact one-form $\lambda$ defines a vector field $R$ on $V$ by 
$R\in\ker d\lambda$ and $\lambda(R)=1$, which 
is called the Reeb vector field. We assume that 
the contact form is Morse in the sense that all closed orbits of the 
Reeb vector field are nondegenerate in the sense of [BEHWZ]; in particular, the set 
of closed Reeb orbits is discrete. The invariants are defined by counting 
$\Ju$-holomorphic curves in $\IR\times V$ which are asymptotically cylindrical over 
chosen collections of Reeb orbits $\Gamma^{\pm}=\{\gamma^{\pm}_1,...,
\gamma^{\pm}_{n^{\pm}}\}$ as the $\IR$-factor tends to $\pm\infty$, see [BEHWZ]. 
The almost complex structure $\Ju$ on the cylindrical 
manifold $\IR\times V$ is required to be cylindrical in the sense that it is  
$\IR$-independent, links the two natural vector fields on $\IR\times V$, namely the 
Reeb vector field $R$ and the $\IR$-direction $\del_s$, by $\Ju\del_s=R$, and turns 
the distribution $\xi$ on $V$ into a complex subbundle of $TV$, 
$\xi=TV\cap \Ju TV$. We denote by $\CM_{g,r}(\Gamma^+,\Gamma^-)$ the corresponding compactified
moduli space of genus $g$ curves with $r$ additional marked points ([BEHWZ],[EGH]). 
Possibly after choosing abstract perturbations using polyfolds (see [HWZ]), obstruction 
bundles ([F2]) or domain-dependent structures ([F1]) following the ideas in [CM] we get that 
$\CM_{g,r}(\Gamma^+,\Gamma^-)$ is a branched-labelled orbifold with boundaries and corners of dimension 
equal to the Fredholm index of the Cauchy-Riemann operator for $\Ju$. 
{\it Note that in the same way as we will not discuss transversality for the general case but just refer to the upcoming papers on polyfolds by Hofer and his co-workers, in what follows we will for simplicity assume that every moduli space is indeed a manifold with boundaries and corners, since we expect that all the upcoming constructions can be generalized in an appropriate way.} \\  
 
Let us now briefly introduce the algebraic formalism of SFT as described in [EGH]: \\
 
Recall that a multiply-covered Reeb orbit $\gamma^k$ is called bad if 
$\CZ(\gamma^k)\neq\CZ(\gamma)\mod 2$, where $\CZ(\gamma)$ denotes the 
Conley-Zehnder index of $\gamma$. Calling a Reeb orbit $\gamma$ {\it good} if it is not bad we assign to every good Reeb orbit $\gamma$ two formal 
graded variables $p_{\gamma},q_{\gamma}$ with grading 
\begin{equation*} 
|p_{\gamma}|=m-3-\CZ(\gamma),|q_{\gamma}|=m-3+\CZ(\gamma) 
\end{equation*} 
when $\dim V = 2m-1$. In order to include higher-dimensional moduli spaces we further assume that a string 
of closed (homogeneous) differential forms $\Theta=(\theta_1,...,\theta_N)$ on $V$ is chosen and assign to 
every $\theta_i\in\Omega^*(V)$ a formal variables $t_i$ 
with grading
\begin{equation*} |t_i|=2 -\deg\theta_i. \end{equation*}  
Finally, let $\hbar$ be another formal variable of degree $|\hbar|=2(m-3)$. \\

Let $\WW$ be the graded Weyl algebra over $\IC$ of power series in the variables 
$\hbar,p_{\gamma}$ and $t_i$ with coefficients which are polynomials in the 
variables $q_{\gamma}$, which is equipped with the associative product $\star$ in 
which all variables super-commute according to their grading except for the 
variables $p_{\gamma}$, $q_{\gamma}$ corresponding to the same Reeb orbit $\gamma$, 
\begin{equation*} [p_{\gamma},q_{\gamma}] = 
                  p_{\gamma}\star q_{\gamma} -(-1)^{|p_{\gamma}||q_{\gamma}|} 
                  q_{\gamma}\star p_{\gamma} = \kappa_{\gamma}\hbar.
\end{equation*}
($\kappa_{\gamma}$ denotes the multiplicity of $\gamma$.) Following [EGH] we further introduce 
the Poisson algebra $\PP$ of formal power series in the variables $p_{\gamma}$ and $t_i$ with   
coefficients which are polynomials in the variables $q_{\gamma}$ with Poisson bracket given by 
\begin{equation*} 
 \{f,g\} = \sum_{\gamma}\kappa_{\gamma}\Bigl(\frac{\del f}{\del p_{\gamma}}\frac{\del g}{\del q_{\gamma}} -
                          (-1)^{|f||g|}\frac{\del g}{\del p_{\gamma}}\frac{\del f}{\del q_{\gamma}}\Bigr).       
\end{equation*}

As in Gromov-Witten theory we want to organize all moduli spaces $\CM_{g,r}(\Gamma^+,\Gamma^-)$
into a generating function $\IH\in\hbar^{-1}\WW$, called {\it Hamiltonian}. In order to include also higher-dimensional 
moduli spaces, in [EGH] the authors follow the approach in Gromov-Witten theory to integrate the chosen differential forms 
$\theta_1,...,\theta_N$ over the moduli spaces after pulling them back under the evaluation map from target manifold $V$. The 
Hamiltonian $\IH$ is then defined by
\begin{equation*}
 \IH = \sum_{\Gamma^+,\Gamma^-} \int_{\CM_{g,r}(\Gamma^+,\Gamma^-)/\IR}
 \ev_1^*\theta_{i_1}\wedge...\wedge\ev_r^*\theta_{i_r}\; \hbar^{g-1}t^Ip^{\Gamma^+}q^{\Gamma^-}
\end{equation*}
with $t^I=t_{i_1}...t_{i_r}$, $p^{\Gamma^+}=p_{\gamma^+_1}...p_{\gamma^+_{n^+}}$ and $q^{\Gamma^-}=q_{\gamma^-_1}...q_{\gamma^-_{n^-}}$. 
Expanding 
\begin{equation*} \IH=\hbar^{-1}\sum_g \IH_g \hbar^g \end{equation*} 
we further get a rational Hamiltonian $\Ih=\IH_0\in\PP$, which counts only curves with genus zero. \\

While the Hamiltonian $\IH$ explicitly depends on the chosen contact form, the cylindrical almost complex structure, the differential forms and abstract polyfold perturbations making all moduli spaces regular, it is outlined in [EGH] how to construct algebraic invariants, which just depend on the contact structure and the cohomology classes of the differential forms. \\
  
\subsection{Gravitational descendants}
For the relation to integrable systems it is outlined in [E] that, as in Gromov-Witten theory, symplectic field theory must be enriched by considering so-called {\it gravitational descendants} of the {\it primary} Hamiltonian $\IH$. \\

Before we give a rigorous definition of gravitational descendants in SFT, we recall the definition from 
Gromov-Witten theory. Denote by $\CM_r=\CM_{g,r}(X,J)$ the compactified moduli space of closed $J$-holomorphic curves in the closed symplectic manifold $X$ of genus $g$ with $r$ marked points (and fixed homology class). Following [MDSa] we introduce over $\CM_r$ so-called {\it tautological line bundles}  $\LL_1,...,\LL_r$, where the fibre of $\LL_i$ over a punctured curve $(u,z_1,...,z_r)\in\CM_r$ in the noncompactified moduli space is given by the cotangent line to the underlying, possibly unstable closed nodal Riemann surface $S$ at the $i$.th marked point,
\begin{equation*} (\LL_i)_{(u,z_1,...,z_r)} = T^*_{z_i}S,\;\;i=1,...,r. \end{equation*}
To be more formal, observe that there exists a canonical map $\pi: \CM_{r+1}\to\CM_r$ by forgetting the $(r+1)$.st marked point and stabilizing the map, where the fibre over the curve $(u,z_1,...,z_r)$ agrees with the curve itself. Then the tautological line bundle $\LL_i$ can be defined as the pull-back of the vertical cotangent line bundle of $\pi: \CM_{r+1}\to\CM_r$ under the canonical section $\sigma_i:\CM_r\to\CM_{r+1}$ mapping to the $i$.th marked point in the fibre. Note that while the vertical cotangent line bundle is rather a sheaf than a true bundle since it becomes singular at the nodes in the fibres, the pull-backs under the canonical sections are indeed true line bundles as the marked points are different from the nodes and hence these sections avoid the singular loci. \\

Denoting by $c_1(\LL_i)$ the first Chern class of the complex line bundle $\LL_i$, one then considers for the descendant potential of Gromov-Witten theory integrals of the form
\begin{equation*} \int_{\CM_r} \ev_1^*\theta_{i_1}\wedge c_1(\LL_1)^{j_1}\wedge ... \wedge \ev_r^*\theta_{i_r}\wedge c_1(\LL_r)^{j_r}, \end{equation*}
where $(i_k,j_k)\in\{1,...,N\}\times\IN$, which can again be organized into a generating function. \\

Like pulling-back cohomology classes from the target manifold, the introduction of the tautological line bundles hence 
has the effect that the generating function also sees the higher-dimensional moduli spaces. On the other hand, 
in contrast to the former, the latter refers to partially fixing the complex structure on the underlying punctured 
Riemann surface. \\

Before we can turn to the definition of gravitational descendants in SFT, it will turn out to be useful to give an alternative definition, where the integration of the powers of the first Chern classes is replaced by considering zero sets of sections. Restricting for notational simplicity to the case with one marked point, we can define by induction over $j\in\IN$ a nested sequence of moduli spaces $\CM^{j+1}_1\subset\CM^j_1\subset\CM_1$ such that 
\begin{equation*} \int_{\CM_1} \ev^*\theta_i \wedge c_1(\LL)^j = \frac{1}{j!} \cdot \int_{\CM^j_1} \ev^*\theta_i. \end{equation*}

For $j=1$ observe that, since the first Chern class of a line bundle agrees with its Euler class, the homology class obtained by integrating $c_1(\LL)$ over the compactified moduli space $\CM_1$ can be represented by the zero set of a generic section $s_1$ in $\LL$. Note that here we use that $\CM_1$ represents a pseudo-cycle and hence has no codimension-one boundary strata. In other words, we find that 
\begin{equation*} \int_{\CM_1} \ev^*\theta_i \wedge c_1(\LL) \;=\; \int_{\CM^1_1} \ev^*\theta_i, \end{equation*}
where $\CM^1_1=s_1^{-1}(0)$. \\

Now consider the restriction of the tautological line bundle $\LL$ to $\CM^{j-1}_1\subset\CM_1$. Instead of describing the integration of powers of the first Chern class in terms of common zero sets of sections in the same line bundle $\LL$, it turns out to be more geometric (see 2.2) to choose a section $s_j$ not in $\LL$ but in its $j$-fold (complex) tensor product $\LL^{\otimes j}$ and define 
\begin{equation*} \CM^j_1 = s_j^{-1}(0) \subset\CM^{j-1}_1. \end{equation*} 
Since $c_1(\LL^{\otimes j}) = j\cdot c_1(\LL)$ it follows that 
\begin{equation*} \int_{\CM^j_1} \ev^*\theta_i = j\cdot \int_{\CM^{j-1}_1} \ev^*\theta_i\wedge c_1(\LL) \end{equation*}
so that by induction
\begin{equation*} \int_{\CM_1} \ev^*\theta_i \wedge c_1(\LL)^j = \frac{1}{j!} \cdot \int_{\CM^j_1} \ev^*\theta_i \end{equation*}
as desired. \\

While the result of the integration is well-known to be independent of the choice of the almost complex structure and the abstract polyfold perturbations, it also follows that the result is independent of the precise choice of the sequence of sections $s_1,...,s_j$.  Like for the almost complex structure and the perturbations this results from the fact that the moduli spaces studied in Gromov-Witten theory have no codimension-one boundary. \\

On the other hand, it is well-known that the moduli spaces in SFT typically have codimension-one boundary, so that now the result of the integration will not only depend on the chosen contact form, cylindrical almost complex structure and abstract polyfold perturbations, but also additionally explicitly depend on the chosen sequences of sections $s_1,...,s_j$. While the Hamiltonian is hence known to depend on all extra choices, it is well-known from Floer theory that we can expect to find algebraic invariants independent of these choices. \\

While the problem of dependency on contact form, cylindrical almost complex structure and abstract polyfold perturbations is sketched in [EGH], we will now show how to include gravitational descendants into their algebraic constructions. For this we will define {\it descendants of moduli spaces}, which we obtain as zero sets of {\it coherent collections of sections} in the tautological line bundles over all moduli spaces. \\

From now on let $\CM_r$ denote the moduli space $\CM_{g,r}(\Gamma^+,\Gamma^-)/\IR$ studied in SFT for chosen collections 
of Reeb orbits $\Gamma^+,\Gamma^-$. In complete analogy to Gromov-Witten theory we can introduce $r$ tautological line 
bundles $\LL_1,...,\LL_r$, where the fibre of $\LL_i$ over a punctured curve $(u,z_1,...,z_r)\in\CM_r$ is again given 
by the cotangent line to the underlying, possibly unstable nodal Riemann surface (without ghost components) at the 
$i$.th marked point and which again formally can be defined as the pull-back of the vertical cotangent line 
bundle of $\pi: \CM_{r+1}\to\CM_r$ under the canonical section $\sigma_i: \CM_r\to\CM_{r+1}$ mapping to the $i$.th marked 
point in the fibre. Note again that while the vertical cotangent line bundle is rather a sheaf than a true bundle since 
it becomes singular at the nodes in the fibres, the pull-backs under the canonical sections are still true line bundles 
as the marked points are different from the nodes and hence these sections avoid the singular loci. \\

For notational simplicity let us again restrict to the case $r=1$. Following the compactness statement in [BEHWZ], the codimension-one boundary of $\CM_1$ consists of curves with two levels (in the sense of [BEHWZ]), whose moduli spaces can be represented as products $\CM_{1,1}\times\CM_{2,0}$ or $\CM_{1,0}\times\CM_{2,1}$ of moduli spaces of strictly lower dimension, where the marked point sits on the first or the second level. As we want to keep the notation as simple as possible, note that here and in what follows for product moduli spaces {\it the first index refers to the level} and not to the genus of the curve. To be more precise, after introducing asymptotic markers as in [EGH] for orientation issues, one obtains a fibre rather than a direct product, see also [F2]. However, since all the bundles and sections we will consider do or should not depend on these asymptotic markers, we will forget about this issue in order to keep the notation as simple as possible. \\ 

On the other hand, it directly follows from the definition of the tautological line bundle $\LL$ over $\CM_1$ that over the boundary components $\CM_{1,1}\times\CM_{2,0}$ and $\CM_{1,0}\times\CM_{2,1}$ it is given by 
\begin{equation*} \LL|_{\CM_{1,1}\times\CM_{2,0}} = \pi_1^*\LL_1,\;\LL|_{\CM_{1,0}\times\CM_{2,1}} = \pi_2^*\LL_2, \end{equation*}
where $\LL_1$, $\LL_2$ denotes the tautological line bundle over the moduli space $\CM_{1,1}$, $\CM_{2,1}$ and $\pi_1$, $\pi_2$ is the projection onto the first or second factor, respectively. \\

With this we can now introduce the notion of coherent collections of sections in (tensor products of) tautological line bundles. \\
\\
{\bf Definition 1.1:} {\it Assume that we have chosen sections $s$ in the tautological line bundles $\LL$ over all moduli spaces $\CM_1$ of $\Ju$-holomorphic curves with one additional marked point. Then this collection of sections $(s)$ is called} coherent {\it if for every section $s$ in $\LL$ over a moduli space $\CM_1$ the following holds: Over every codimension-one boundary component $\CM_{1,1}\times\CM_{2,0}$, $\CM_{1,0}\times\CM_{2,1}$ of $\CM_1$ the section $s$ agrees with the pull-back $\pi_1^*s_1$, $\pi_2^*s_2$ of the chosen section $s_1$, $s_2$ in the tautological line bundle $\LL_1$ over $\CM_{1,1}$, $\LL_2$ over $\CM_{2,1}$, respectively.}  \\

\noindent{\bf Remark:} Since in the end we will again be interested in the zero sets of these sections, {\it we will assume that all occuring sections are transversal to the zero section.} Furthermore, we want to assume that all the chosen sections are indeed {\it invariant under the obvious symmetries like reordering of punctures and marked points.} In order to meet both requirements, it follows that actually need to employ {\it multi-sections} as in [CMS], which we however want to suppress for the rest of this exposition. \\

The important observation is clearly that one can always find coherent collections of (transversal) sections $(s)$ by using induction on the dimension of the underlying moduli space. While for the induction start it suffices to choose a non-vanishing section in the tautological line bundle over the moduli space of orbit cylinders with one marked point, for the induction step observe that the coherency condition fixes the section on the boundary of the moduli space. Here it is important to remark that the coherency condition further ensures that two different codimension-one boundary components actually agree on their common boundary strata of higher codimension. On the other hand, we can use our assumption that every moduli space is indeed a manifold with corners to obtain the desired section by simply extending the section from the boundary to the interior of the moduli space in an arbitrary way. \\
  
For a given coherent collection of transversal sections $(s)$ we will again define for every moduli space 
\begin{equation*} \CM^1_1 = s^{-1}(0) \subset \CM_1. \end{equation*}
As an immediate consequence of the above definition we find that $\CM^1_1$ is a neat submanifold (with corners) of $\CM_1$, i.e., the components of the codimension-one boundary of $\CM^1_1$ are given by products $\CM^1_{1,1}\times\CM_{2,0}$ and $\CM_{1,0}\times\CM^1_{2,1}$, where $\CM^1_{1,1} = s_1^{-1}(0)$, $\CM^1_{2,1} = s_2^{-1}(0)$  for the section $s_1$ in $\LL_1$ over $\CM_{1,1}$, $s_2$ in $\LL_2$ over $\CM_{2,1}$, respectively. To be more precise, since we actually need to work with multi-sections rather than sections in the usual sense, the zero set is indeed a branched-labelled manifold. On the other hand, since we already suppressed the fact that our moduli spaces are indeed branched and labelled, we want to continue ignoring this technical aspect. On the other hand, we can use the above result as an induction start to obtain for every moduli space $\CM_1$ a sequence of nested subspaces $\CM^j_1\subset\CM^{j-1}_1\subset\CM_1$ as in Gromov-Witten theory. \\
\\
{\bf Definition 1.2:} {\it Let $j\in\IN$. Assume that for all moduli spaces we have chosen $\CM^{j-1}_1\subset\CM_1$ such that the components of the codimension-one boundary of $\CM^{j-1}_1$ are given by products of the form $\CM^{j-1}_{1,1}\times\CM_{2,0}$ and $\CM_{1,0}\times\CM^{j-1}_{2,0}$. Then we again call a collection of transversal sections $(s_j)$ in the j-fold tensor products $\LL^{\otimes j}$ of the tautological line bundles over $\CM^{j-1}_1\subset \CM_1$} coherent {\it if for every section $s_j$ the following holds: Over every codimension-one boundary component $\CM^{j-1}_{1,1}\times\CM_{2,0}$, $\CM_{1,0}\times\CM^{j-1}_{2,1}$ of $\CM^{j-1}_1$ the section $s_j$ agrees with the pull-back $\pi_1^*s_{1,j}$, $\pi_2^*s_{2,j}$ of the section $s_{1,j}$, $s_{2,j}$ in the line bundle $\LL^{\otimes j}_1$ over $\CM^{j-1}_{1,1}$, $\LL^{\otimes j}_2$ over $\CM^{j-1}_{2,1}$, respectively.} \\

With this we will now introduce {\it (gravitational) descendants of moduli spaces.} \\
\\
{\bf Definition 1.3:} {\it Assume that we have inductively defined a subsequence of nested subspaces $\CM^j_1\subset\CM^{j-1}_1\subset\CM_1$ by requiring that 
$\CM^j_1 = s_j^{-1}(0)\subset\CM^{j-1}_1$ for a coherent collection of sections $s_j$ in the line bundles $\LL^{\otimes j}$ over the moduli spaces $\CM^{j-1}_1$. Then we call $\CM^j_1$ the $j$.th (gravitational) descendant of $\CM_1$.} \\

Let $\WW^0$ be the graded Weyl algebra over $\IC$ of power series in the variables $\hbar$ and $p_{\gamma}$ with coefficients which are polynomials in the 
variables $q_{\gamma}$, which is obtained from the big Weyl algebra $\WW$ by setting all variables $t_i$ equal to zero. In the same way define the subalgebra $\PP^0$ of the Poisson algebra $\PP$. Apart from the Hamiltonian $\IH^0\in\hbar^{-1}\WW^0$ counting only curves with no additional marked points,
\begin{equation*} \IH^0 = \sum_{\Gamma^+,\Gamma^-} \#\CM_{g,0}(\Gamma^+,\Gamma^-)/\IR\; \hbar^{g-1} p^{\Gamma^+}q^{\Gamma^-}, \end{equation*}
we now want to use the chosen differential forms $\theta_i\in\Omega^*(V)$, $i=1,...,N$ and the sequences $\CM^j_1 = \CM^j_{g,1}(\Gamma^+,\Gamma^-)/\IR$ of gravitational descendants to define sequences of new SFT Hamiltonians $\IH^1_{i,j}\in\hbar^{-1}\WW^0$, $(i,j)\in\{1,...,N\}\times \IN$, by
\begin{equation*} \IH^1_{i,j} = \sum_{\Gamma^+,\Gamma^-} \int_{\CM^j_{g,1}(\Gamma^+,\Gamma^-)/\IR}\ev^*\theta_i\; \hbar^{g-1} p^{\Gamma^+}q^{\Gamma^-}. \end{equation*}

We want to emphasize that the following statement is not yet a theorem in the strict mathematical sense as the analytical foundations of symplectic field theory, in particular, the neccessary transversality theorems for the Cauchy-Riemann operator, are not yet fully established. Since it can be expected that the polyfold project by Hofer and his collaborators sketched in [HWZ] will provide the required transversality theorems, we follow other papers in the field in proving everything up to transversality and state it nevertheless as a theorem. \\
\\
{\bf Theorem 1.4:} {\it Counting holomorphic curves with one marked point after integrating differential forms and introducing gravitational descendants defines a sequence of distinguished elements} 
\begin{equation*} \IH^1_{i,j}\in H_*(\hbar^{-1}\WW^0,D^0) \end{equation*} 
{\it in the full SFT homology algebra with differential $D^0=[\IH^0,\cdot]: \hbar^{-1}\WW^0\to\hbar^{-1}\WW^0$, which commute with respect to the bracket on $H_*(\hbar^{-1}\WW^0,D^0)$,} 
\begin{equation*} [\IH^1_{i,j},\IH^1_{k,\ell}] = 0,\; (i,j),(k,\ell)\in\{1,...,N\}\times\IN. \end{equation*}
\\
{\it Proof:} While the boundary equation $D^0\circ D^0=0$ is well-known to follow from the identity $[\IH^0,\IH^0]=0$, the fact that every $\IH^1_{i,j}$, $(i,j)\in\{1,...,N\}\times\IN$ defines an element in the homology $H_*(\hbar^{-1}\WW^0,D^0)$ follows from the identity
\begin{equation*} [\IH^0,\IH^1_{i,j}]=0,\end{equation*} 
since this proves $\IH^1_{i,j}\in\ker D^0$. On the other hand, in order to see that any two $\IH^1_{i,j}$, $\IH^1_{k,\ell}$ commute {\it after passing to homology} it suffices to prove the identity
\begin{equation*} [\IH^1_{i,j},\IH^1_{k,\ell}]\pm[\IH^0,\IH^2_{(i,j),(k,\ell)}]=0 \end{equation*} 
for any $(i,j),(k,\ell)\in\{1,...,N\}\times\IN$, where the new Hamiltonian $\IH^2_{(i,j),(k,\ell)}$ is defined below using descendant moduli spaces with two additional marked points. \\

The latter two identities directly follow from our definition of gravitational descendants of moduli spaces based on the definition of coherent sections in tautological line bundles and the compactness theorem in [BEHWZ]. Indeed, in the same way as the identity $[\IH^0,\IH^0]=0$ follows from the fact that the codimension-one boundary of every moduli space $\CM_0$ is formed by products of moduli spaces $\CM_{1,0}\times\CM_{2,0}$, the second identity $[\IH^0,\IH^1_{i,j}]=0$ follows from the fact that the codimension-one boundary of a descendant moduli space $\CM^j_1$ is given by products of the form $\CM^j_{1,1}\times\CM_{2,0}$ and $\CM_{1,0}\times\CM^j_{2,1}$. \\

In order to prove the third identity $[\IH^1_{i,j},\IH^1_{k,\ell}]\pm[\IH^0,\IH^2_{(i,j),(k,\ell)}]=0$ for every $(i,j),(k,\ell)\in\{1,...,N\}\times\IN$, we slightly have to enlarge our definition of gravitational descendants in order to include moduli spaces with two additional marked points. For this observe that for every pair $j,k\in\IN$ we can define decendants $\CM^{(j,k)}_2$ of $\CM_2$ by setting $\CM^{(j,k)}_2 = \CM^{(j,0)}_2\cap\CM^{(0,k)}_2$, where $\CM^{(j,0)}_2, \CM^{(0,k)}_2\subset\CM_2$ are defined in the same way as $\CM^j_1,\CM^k_1\subset\CM_1$ by simply forgetting the second or first additional marked point, respectively. Since the boundary of a moduli space of curves with two marked points consists of products of the form $\CM_{1,1}\times\CM_{2,1}$ and $\CM_{1,0}\times\CM_{2,2}$, $\CM_{1,2}\times\CM_{2,0}$, it follows that the boundary of $\CM^{(j,0)}_2$ consists of products $\CM^j_{1,1}\times\CM_{2,1}$, $\CM_{1,1}\times\CM^j_{2,1}$ and $\CM_{1,0}\times\CM^{(j,0)}_{2,2}$, $\CM^{(j,0)}_{1,2}\times\CM_{2,0}$. Together with the similar result about the boundary of $\CM^{(0,k)}_2$ and using the inclusions we hence obtain that the codimension-one boundary of $\CM^{(j,k)}_2$ is given by products of the form $\CM^j_{1,1}\times\CM^k_{2,1}$, $\CM^k_{1,1}\times\CM^j_{2,1}$ and $\CM_{1,0}\times\CM^{(j,k)}_{2,2}$, $\CM^{(j,k)}_{1,2}\times\CM_{2,0}$. While summing over the first two products (with signs) we obtain $[\IH^1_{i,j},\IH^1_{k,\ell}]$, summing over the latter two we get $[\IH^0,\IH^2_{(i,j),(k,\ell)}]$, which hence sum up to zero. $\qed$ \\  
\\
{\bf Remark:} While the proof suggests that for the above algebraic relations one only has to care about the codimension-one boundary strata of the moduli spaces, it is actually even more important that the coherency condition further ensures that two different codimension-one boundary components can be glued along their common boundary strata of higher codimension. \\

As above we further again obtain a rational version of the above statement by expanding $\IH^0$ and the $\IH^1_{i,j}$ in powers 
of $\hbar$. \\
\\
{\bf Corollary 1.5:} {\it Counting rational holomorphic curves with one marked point after integrating differential forms and introducing gravitational descendants defines a sequence of distinguished elements} 
\begin{equation*} \Ih^1_{i,j}\in H_*(\PP^0,d^0), \end{equation*} 
{\it in the rational SFT homology algebra with differential $d^0=\{\Ih^0,\cdot\}: \PP^0\to\PP^0$, which commute with respect to the Poisson bracket on $H_*(\PP^0,d^0)$,} 
\begin{equation*} \{\Ih^1_{i,j},\Ih^1_{k,\ell}\} = 0,\; (i,j),(k,\ell)\in\{1,...,N\}\times\IN. \end{equation*}

So far we have only considered the case with one additional marked point. On the other hand, the general case with $r$ additional marked points is just notationally more involved. Indeed, as we did in the proof of the above theorem we can easily define for every moduli space $\CM_r$ with $r$ additional marked points and every $r$-tuple of natural numbers $(j_1,...,j_r)$ descendants $\CM^{(j_1,...,j_r)}_r\subset\CM_r$ by setting
\begin{equation*} \CM^{(j_1,...,j_r)}_r = \CM^{(j_1,0,...,0)}_r\cap ... \cap \CM^{(0,...,0,j_r)}_r, \end{equation*}
where the descendant moduli spaces $\CM^{(0,...,0,j_k,0,...,0)}_r\subset\CM_r$ are defined in the same way as the one-point descendant $\CM^{j_k}_1\subset\CM_1$ by looking at the $r$ tautological line bundles over the moduli space $\CM_r = \CM_r(\Gamma^+,\Gamma^-)/\IR$ separately and forgetting about the other points. \\

With this we can define the descendant Hamiltonian of SFT, which we will continue denoting by $\IH$, while the Hamiltonian defined in [EGH] will from now on be called {\it primary}. In order to keep track of the descendants we will assign to every chosen differential form $\theta_i$ now a sequence of formal variables $t_{i,j}$ with grading 
\begin{equation*} |t_{i,j}|=2(1-j) -\deg\theta_i. \end{equation*} 
Then the descendant Hamiltonian $\IH$ of SFT is defined by 
\begin{equation*}
 \IH = \sum_{\Gamma^+,\Gamma^-,I} \int_{\CM^{(j_1,...,j_r)}_{g,r}(\Gamma^+,\Gamma^-)/\IR}
 \ev_1^*\theta_{i_1}\wedge...\wedge\ev_r^*\theta_{i_r}\; \hbar^{g-1}t^Ip^{\Gamma^+}q^{\Gamma^-},
\end{equation*}
where $p^{\Gamma^+}=p_{\gamma^+_1} ... p_{\gamma^+_{n^+}}$, $q^{\Gamma^-}=q_{\gamma^-_1} ... q_{\gamma^-_{n^-}}$ and $t^I=t_{i_1,j_1} ... t_{i_r,j_r}$ for $I=((i_1,j_1),...,(i_r,j_r))$. \\

Note that expanding the Hamiltonian $\IH$ in powers of the formal variables $t_{i,j}$,
\begin{equation*} 
  \IH = \IH^0 + \sum_{i,j} t_{i,j}\IH^1_{i,j} + o(t^2), 
\end{equation*}
we get back our Hamiltonians $\IH^0$ and the sequences of descendant Hamiltonians $\IH^1_{i,j}$ from above and it is easy to see that the primary Hamiltonian from [EGH] is recovered by setting all formal variables $t_{i,j}$ with $j>0$ equal to zero. \\

In the same way as it was shown for the primary Hamiltonian in [EGH], the descendant Hamiltonian continues to satisfy the master equation $[\IH,\IH]=0$, which is just a generalization of the identities for $\IH^0$, $\IH^1_{i,j}$ and hence can be shown along the same lines by studying the codimension-one boundaries of descendant moduli spaces. On the other hand, expanding $\IH\in\hbar^{-1}\WW$ in terms of powers of $\hbar$, 
\begin{equation*} \IH=\sum_g \hbar^{g-1}\IH_g,  \end{equation*} 
note that for the rational descendant Hamiltonian $\Ih=\IH_0\in\PP$ we still have $\{\Ih,\Ih\}=0$. \\

\subsection{Invariance statement}

We now turn to the question of independence of these nice algebraic structures from the choices like contact form, cylindrical almost complex structure, abstract polyfold perturbations and, of course, the choice of the coherent collection of sections. This is the content of the following theorem, where we however again want to emphasize that the following statement is not yet a theorem in the strict mathematical sense as the analytical foundations of symplectic field theory, in particular, the neccessary transversality theorems for the Cauchy-Riemann operator, are not yet fully established. \\   
\\
{\bf Theorem 1.6:} {\it For different choices of contact form $\lambda^{\pm}$, cylindrical almost complex structure $\Ju^{\pm}$ , abstract polyfold perturbations and sequences of coherent collections of sections $(s^{\pm}_j)$ the resulting systems of commuting operators $\IH^{1,-}_{i,j}$ on $H_*(\hbar^{-1}\WW^{0,-},D^{0,-})$ and $\IH^{1,+}_{i,j}$ on $H_*(\hbar^{-1}\WW^{0,+},D^{0,+})$ are isomorphic, i.e., there exists an isomorphism of the Weyl algebras $H_*(\hbar^{-1}\WW^{0,-},D^{0,-})$ and $H_*(\hbar^{-1}\WW^{0,+},D^{0,+})$ which maps $\IH^{1,-}_{i,j}\in H_*(\hbar^{-1}\WW^{0,-},D^{0,-})$ to $\IH^{1,+}_{i,j}\in H_*(\hbar^{-1}\WW^{0,+},D^{0,+})$.} \\

As above we clearly also get a rational version of the invariance statement: \\
\\
{\bf Corollary 1.7:} {\it For different choices of contact form $\lambda^{\pm}$, cylindrical almost complex structure $\Ju^{\pm}$, abstract polyfold perturbations and sequences of coherent collections of sections $(s^{\pm}_j)$ the resulting system of Poisson-commuting functions $\Ih^{1,-}_{i,j}$ on $H_*(\PP^{0,-},d^{0,-})$ and $\Ih^{1,+}_{i,j}$ on $H_*(\PP^{0,+},d^{0,+})$ are isomorphic, i.e., there exists an isomorphism of the Poisson algebras $H_*(\PP^{0,-},d^{0,-})$ and $H_*(\PP^{0,+},d^{0,+})$ which maps $\Ih^{1,-}_{i,j}\in H_*(\PP^{0,-},d^{0,-})$ to $\Ih^{1,+}_{i,j}\in H_*(\PP^{0,+},d^{0,+})$.} \\

This theorem is an extension of the theorem in [EGH] which states that for different choices of auxiliary data the small Weyl algebras $H_*(\hbar^{-1}\WW^{0,-},D^{0,-})$ and $H_*(\hbar^{-1}\WW^{0,+},D^{0,+})$ are isomorphic. On the other hand, assuming that the contact form, the cylindrical almost complex structure and also the abstract polyfold sections are fixed to have well-defined moduli spaces, the isomorphism of the homology algebras is the identity and hence the theorem states the sequence of commuting operators is indeed independent of the chosen sequences of coherent collections of sections $(s^{\pm}_j)$,  
\begin{equation*} \IH^{1,-}_{i,j} = \IH^{1,+}_{i,j}\in H_*(\hbar^{-1}\WW^0,D^0). \end{equation*}
$ $\\

For the proof we have to extend the proof in [EGH] to include gravitational descendants. To this end we have to study sections in the tautological line bundles over moduli spaces of holomorphic curves in symplectic manifolds with cylindrical ends. \\
 
Let $(W,\omega)$ be a symplectic manifold with cylindrical ends $(\IR^+\times V^+,\lambda^+)$ and 
$(\IR^-\times V^-,\lambda^-)$ in the sense of [BEHWZ] which is equipped with an almost complex structure 
$\Ju$ which agrees with the cylindrical almost complex structures $\Ju^{\pm}$ on $\IR^+\times V^+$. Then we 
study $\Ju$-holomorphic curves in $W$ which are asymptotically cylindrical over 
chosen collections of orbits $\Gamma^{\pm}=\{\gamma^{\pm}_1,...,
\gamma^{\pm}_{n^{\pm}}\}$ of the Reeb vector fields $R^{\pm}$ in $V^{\pm}$ as the $\IR^{\pm}$-factor tends 
to $\pm\infty$, see [BEHWZ], and denote by $\IM_{g,r}(\Gamma^+,\Gamma^-)$ the corresponding 
moduli space of genus $g$ curves with $r$ additional marked points ([BEHWZ],[EGH]). 
Possibly after choosing abstract perturbations using polyfolds, obstruction 
bundles or domain-dependent structures, which agree with chosen abstract perturbations in the boundary as described above, we find that 
$\IM_{g,r}(\Gamma^+,\Gamma^-)$ is a weighted branched manifold of dimension equal to the Fredholm index of the Cauchy-Riemann operator for $\Ju$. Note that as remarked above we will for simplicity assume that moduli space is indeed a manifold with corners, since this will be sufficient for our example and we expect that all the upcoming 
constructions can be generalized in an appropriate way. We further extend the chosen differential forms $\theta^{\pm}_1,...\theta^{\pm}_N$ on $V^{\pm}$ to differential forms $\theta_1,...,\theta_N$ on $W$ as described in [EGH]. \\

From now on let $\CM_r$ denote the moduli space $\CM_{g,r}(\Gamma^+,\Gamma^-)$ of holomorphic curves in $W$ for chosen collections of Reeb orbits $\Gamma^+,\Gamma^-$. Note in particular that there is no longer an $\IR$-action on the moduli space which we have to quotient out. In order to distinguish these moduli spaces in non-cylindrical manifolds from those of holomorphic curves in the cylindrical manifolds, we will use the short-hand notation $\CM^{\pm}_r$ for moduli spaces $\CM_{g,r}(\Gamma^+,\Gamma^-)/\IR$ of holomorphic curves in $\IR\times V^{\pm}$, respectively. Like in Gromov-Witten theory we can introduce $r$ tautological line bundles $\LL_1,...,\LL_r$, where the fibre of $\LL_i$ over a punctured curve $(u,z_1,...,z_r)\in\IM_r$ in the noncompactified moduli space is again given by the cotangent line to the underlying closed Riemann surface at the $i$.th marked point and which formally can be defined as the pull-back of the vertical cotangent line bundle under the canonical section $\sigma_i$ of $\pi: \CM_{r+1}\to\CM_r$ mapping to the $i$.th marked point in the fibre. \\

For notational simplicity let us again restrict to the case $r=1$. Following the compactness statement in [BEHWZ] the codimension-one boundary of $\CM_1$ now consists of curves with one non-cylindrical level and one cylindrical level (in the sense of [BEHWZ]), whose moduli spaces can now be represented as products $\CM_{1,1}\times\CM^+_{2,0}$, $\CM^-_{1,1}\times\CM_{2,0}$ or $\CM_{1,0}\times\CM^+_{2,1}$, $\CM^-_{1,0}\times\CM_{2,1}$ of moduli spaces of strictly lower dimension, where the marked point sits on the first or the second level. Again note that here and in what follows for product moduli spaces the first index refers to the level and not to the genus of the curve. Furthermore it follows from the definition of the tautological line bundle $\LL$ over $\CM_1$ that over the boundary components $\CM_{1,1}\times\CM^+_{2,0}$, $\CM^-_{1,1}\times\CM_{2,0}$ and $\CM_{1,0}\times\CM^+_{2,1}$, $\CM^-_{1,0}\times\CM_{2,1}$ it is given by 
\begin{eqnarray*} 
  \LL|_{\CM_{1,1}\times\CM^+_{2,0}} = \pi_1^*\LL_1,&& \LL|_{\CM_{1,0}\times\CM^+_{2,1}} = \pi_2^*\LL^+_2, \\
  \LL|_{\CM^-_{1,1}\times\CM_{2,0}} = \pi_1^*\LL^-_1,&&\LL|_{\CM^-_{1,0}\times\CM_{2,1}} = \pi_2^*\LL_2, 
\end{eqnarray*}
where $\LL^{(-)}_1$, $\LL^{(+)}_2$ denotes the tautological line bundle over the moduli space $\CM^{(-)}_{1,1}$, $\CM^{(+)}_{2,1}$ and $\pi_1$, $\pi_2$ is the projection onto the first or second factor, respectively. \\

With this we can now introduce collections of sections in (tensor products of) tautological line bundles coherently connecting two chosen coherent collections of sections. \\
\\
{\bf Definition 1.8:} {\it Let $W$ be a symplectic manifold with cylindrical ends $V^{\pm}$ and let $(s_{\pm})$ be two coherent collections of sections in the tautological line bundles $\LL^{\pm}$ over all moduli spaces $\CM^{\pm}_1$ of $\Ju$-holomorphic curves with one additional marked point in the cylindrical manifolds $\IR\times V^{\pm}$. Assume that we have chosen transversal sections $s$ in the tautological line bundles $\LL$ over all moduli spaces $\CM_1$ of $\Ju$-holomorphic curves in the non-cylindrical manifold $W$ with one additional marked point. Then this collection of sections $(s)$ is called} coherently connecting $(s_-)$ and $(s_+)$ {\it if for every section $s$ in $\LL$ over a moduli space $\CM_1$ the following holds: Over every codimension-one boundary component $\CM_{1,1}\times\CM^+_{2,0}$, $\CM^-_{1,1}\times\CM_{2,0}$ and $\CM_{1,0}\times\CM^+_{2,1}$, $\CM^-_{1,0}\times\CM_{2,1}$ of $\CM_1$ the section $s$ agrees with the pull-back $\pi_1^*s_1$, $\pi_1^*s^-_1$ or $\pi_2^*s^+_2$, $\pi_2^*s_2$ of the chosen sections $s_{1,(-)}$, $s_{2,(+)}$ in the tautological line bundles $\LL^{(-)}_1$ over $\CM^{(-)}_{1,1}$, $\LL^{(+)}_2$ over $\CM^{(+)}_{2,1}$, respectively.}  \\

Note that one can always find collections of sections $(s)$ coherently connecting given coherent collections of sections $(s_+)$ and $(s_-)$ as before by using induction on the dimension of the underlying moduli space. Indeed, for the induction step observe that the coherency condition again fixes the section on the boundary of the moduli space, so that the desired section can be obtained by simply extending the section from the boundary to the interior of the moduli space in an arbitrary way. \\
  
For a given coherently connecting collection of sections $(s)$ we will again define for every moduli space 
\begin{equation*} \CM^1_1 = s^{-1}(0) \subset \CM_1. \end{equation*}
As an immediate consequence of the above definition we find that the components of the codimension-one boundary of $\CM^1_1$ are given by products 
$\CM^1_{1,1}\times\CM^+_{2,0}$, $\CM^{1,-}_{1,1}\times\CM_{2,0}$ and $\CM_{1,0}\times\CM^{1,+}_{2,1}$, $\CM^-_{1,0}\times\CM^1_{2,1}$, where $\CM^{1,(-)}_{1,1} = s_{1,(-)}^{-1}(0)$, $\CM^{1,(+)}_{2,1} = s_{2,(+)}^{-1}(0)$  for the section $s_{1,(-)}$ in $\LL^{(-)}_1$ over $\CM^{(-)}_{1,1}$, $s_{2,(+)}$ in $\LL^{(+)}_2$ over $\CM^{(+)}_{2,1}$, respectively. As before we can use this result as an induction start to obtain for every moduli space $\CM_1$ a sequence of nested subspaces $\CM^j_1\subset\CM^{j-1}_1\subset\CM_1$. \\
\\
{\bf Definition 1.9:} {\it Let $j\in\IN$ and let $(s_{j,\pm})$ be two coherent collections of sections in the $j$-fold tensor products $\LL^{\pm,\otimes j}$ of the tautological line bundles over the $j-1$.st gravitational descendants $\CM^{j-1,\pm}_1 \subset\CM^{\pm}_1$ of all moduli spaces of curves in the cylindrical manifolds $\IR\times V^{\pm}$. Assume that for all moduli spaces of curves in the non-cylindrical manifold $W$ we have chosen $\CM^{j-1}_1\subset\CM_1$ such that the components of the codimension-one boundary of $\CM^{j-1}_1$ are given by products of the form $\CM^{j-1}_{1,1}\times\CM^+_{2,0}$, $\CM^{j-1,-}_{1,1}\times\CM_{2,0}$ and $\CM_{1,0}\times\CM^{j-1,+}_{2,0}$, $\CM^-_{1,0}\times\CM^{j-1}_{2,0}$. Then we again call a collection of transversal sections $(s_j)$ in the j-fold tensor products $\LL^{\otimes j}$ of the tautological line bundles over $\CM^{j-1}_1\subset \CM_1$} coherently connecting $(s_{j,-})$ and $(s_{j,+})$ {\it if for every section $s_j$ the following holds: Over every codimension-one boundary component $\CM^{j-1}_{1,1}\times\CM^+_{2,0}$, $\CM^{j-1,-}_{1,1}\times\CM_{2,0}$ and $\CM_{1,0}\times\CM^{j-1,+}_{2,1}$, $\CM^-_{1,0}\times\CM^{j-1}_{2,1}$ of $\CM^{j-1}_1$ the section $s_j$ agrees with the pull-back $\pi_1^*s_{1,j}$, $\pi_1^*s_{1,j,-}$ or $\pi_2^*s_{2,j,+}$, $\pi_2^*s_{2,j}$ of the section $s_{1,j,(-)}$, $s_{2,j,(+)}$ in the line bundle $\LL^{(-),\otimes j}_1$ over $\CM^{j-1,(-)}_{1,1}$, $\LL^{(+),\otimes j}_2$ over $\CM^{j-1,(+)}_{2,1}$, respectively.} \\

With this we can now introduce gravitational descendants of moduli spaces for symplectic manifolds with cylindrical ends. \\
\\
{\bf Definition 1.10:} {\it Assume that we have inductively defined subsequence of nested subspaces $\CM^j_1\subset\CM^{j-1}_1\subset\CM_1$ by requiring that 
$\CM^j_1 = s_j^{-1}(0)\subset\CM^{j-1}_1$ for a collection of sections $s_j$ in the line bundles $\LL^{\otimes j}$ over the moduli spaces $\CM^{j-1}_1$ coherently connecting the coherent collections of sections $(s_{j,-})$ and $(s_{j,+})$. Then we call $\CM^j_1$ the $j$.th (gravitational) descendant of $\CM_1$.} \\
 
In order to prove the above invariance theorem we now recall the extension of the algebraic formalism of SFT from cylindrical manifolds to symplectic cobordisms with cylindrical ends as described in [EGH]. \\

Let $\DD^0$ be the space of formal power series in the variables $\hbar,p^+_{\gamma}$ with coefficients which are polynomials in the 
variables $q^-_{\gamma}$. Elements in $\WW^{0,\pm}$ then act as differential operators from the right/left on 
$\DD^0$ via the replacements 
\begin{equation*} 
  q^+_{\gamma}\mapsto \kappa_{\gamma}\hbar\overleftarrow{\frac{\del}{\del p^+_{\gamma}}},\;\;
  p^-_{\gamma}\mapsto\kappa_{\gamma}\hbar\overrightarrow{\frac{\del}{\del q^-_{\gamma}}}.
\end{equation*}

Apart from the potential $\IF^0\in\hbar^{-1}\WW^0$ counting only curves in $W$ with no additional marked points,
\begin{equation*} \IF^0 = \sum_{\Gamma^+,\Gamma^-} \#\CM_{g,0}(\Gamma^+,\Gamma^-)\; \hbar^{g-1} p^{\Gamma^+}q^{\Gamma^-}, \end{equation*}
we now want to use the extensions $\theta_i$, $i=1,...,N$ on $W$ of the chosen differential forms $\theta^{\pm}_1,...\theta^{\pm}_N$ on $V^{\pm}$ and these sequences 
$\CM^j_1 = \CM^j_{g,1}(\Gamma^+,\Gamma^-)$ of gravitational descendants to define sequences of new SFT potentials $\IF^1_{i,j}$, $(i,j)\in\{1,...,N\}\times \IN$, by
\begin{equation*} \IF^1_{i,j} = \sum_{\Gamma^+,\Gamma^-} \int_{\CM^j_{g,1}(\Gamma^+,\Gamma^-)}\ev^*\theta_i\; \hbar^{g-1} p^{\Gamma^+}q^{\Gamma^-}. \end{equation*}

For the potential counting curves with no additional marked points we have the following identity, where we however again want to emphasize that the following statement should again be understood as a theorem up to the transversality problem in SFT. \\
$ $\\
{\bf Theorem ([EGH]):} {\it The potential $\IF^0\in\hbar^{-1}\DD$ satisfies the master equation}
\begin{equation*} e^{\IF^0}\overleftarrow{\IH^{0,+}} - \overrightarrow{\IH^{0,-}}e^{\IF^0} = 0. \end{equation*}

In [EGH] it is shown that this implies that 
\begin{equation*}
  D^{\IF^0}: \hbar^{-1}\DD^0\to\hbar^{-1}\DD^0,\; 
  D^{\IF^0}g = e^{-\IF^0}\overrightarrow{\IH^{0,-}}(ge^{\IF^0}) - (-1)^{|g|}(ge^{\IF^0})\overleftarrow{\IH^{0,+}}e^{-\IF^0} 
\end{equation*}
satisfies $D^{\IF^0}\circ D^{\IF^0} = 0$ and hence can be used to define the homology algebra 
$H_*(\hbar^{-1}\DD^0,D^{\IF^0})$. Furthermore it is shown that the maps 
\begin{eqnarray*}
  &&F^{0,-}: \hbar^{-1}\WW^{0,-}\to\hbar^{-1}\DD^0,\; f\mapsto e^{-\IF^0}\overrightarrow{f}e^{+\IF^0}, \\
  &&F^{0,+}: \hbar^{-1}\WW^{0,+}\to\hbar^{-1}\DD^0,\; f\mapsto e^{+\IF^0}\overleftarrow{f}e^{-\IF^0}
\end{eqnarray*}
commute with the boundary operators,
\begin{equation*} F^{0,\pm}\circ D^{0,\pm} = D^{\IF^0}\circ F^{0,\pm}, \end{equation*} 
and hence descend to maps between the homology algebras
\begin{equation*} F^{0,\pm}_*: H_*(\hbar^{-1}\WW^{0,\pm},D^{0,\pm})\to H_*(\hbar^{-1}\DD^0,D^{\IF^0}). \end{equation*}

Now assume that the contact forms $\lambda^+$ and $\lambda^-$ are chosen such that they define the same contact structure $(V^+,\xi^+)=(V^-,\xi^-)=:(V,\xi)$ and let $W=\IR\times V$ be the topologically trivial cobordism. Then in [EGH] the authors prove (up to transversality) the following fundamental theorem. \\
\\
{\bf Theorem ([EGH]):} {\it The map}
\begin{equation*} (F^{0,+}_*)^{-1}\circ F^{0,-}_*: H_*(\hbar^{-1}\WW^{0,-},D^{0,-})\to H_*(\hbar^{-1}\WW^{0,+},D^{0,+}) \end{equation*}
{\it is an isomorphism of graded Weyl algebras.}\\
  
For the proof of the invariance statement we want to show that this map identifies the sequences $\IH^{1,\pm}_{i,j}$, $(i,j)\in\{1,...,N\}\times\IN$ on $H_*(\hbar^{-1}\WW^{0,\pm},D^{0,\pm})$. In order to get the right idea for the proof, it turns out to be useful to even enlarge the picture as follows. \\

Precisely in the same way as for cylindrical manifolds we can define for every tuple $(j_1,...,j_r)$ of natural numbers gravitational descendants $\CM^{(j_1,...,j_r)}\subset \CM_1$ of moduli spaces of curves in non-cylindrical manifolds with more than one additional marked point, which are collected in the {\it descendant potential} $\IF\in\hbar^{-1}\DD$, where $\DD$ is again obtained from $\DD^0$ by considering coefficients which are formal powers in the graded formal variables $t_{i,j}$, $(i,j)\in\{1,...,N\}\times\IN$. \\

Assuming for the moment that we have proven the fundamental identity  
\begin{equation*} e^{\IF}\overleftarrow{\IH^+} - \overrightarrow{\IH^-}e^{\IF} = 0 \end{equation*}
and expanding the potential $\IF\in\hbar^{-1}\DD$ and the two Hamiltonians $\IH^{\pm}\in\hbar^{-1}\WW^{\pm}$ in 
powers of the $t$-variables,
\begin{equation*} 
  \IF = \IF^0 + \sum_{i,j} t_{i,j}\IF^1_{i,j} + o(t^2), \;\;
  \IH^{\pm} = \IH^{0,\pm} + \sum_{i,j} t_{i,j}\IH^{1,\pm}_{i,j} + o(t^2), \\
\end{equation*}
we can deduce besides the master equation for $\IF^0$,
\begin{equation*} 
   e^{\IF^0}\overleftarrow{\IH^{0,+}} - \overrightarrow{\IH^{0,-}}e^{\IF^0} = 0 
\end{equation*}
and other identities also the identity  
\begin{equation*} 
   e^{\IF^0}\overleftarrow{\IH^{1,+}_{i,j}} - \overrightarrow{\IH^{1,-}_{i,j}}e^{\IF^0} = 
   \overrightarrow{\IH^{0,-}}(e^{\IF^0}\IF^1_{i,j}) - (e^{\IF^0}\IF^1_{i,j})\overleftarrow{\IH^{0,+}}, 
\end{equation*} 
about $\IF^0$, $\IF^1_{i,j}$ and $\IH^{0,\pm}$, $\IH^{1,\pm}_{i,j}$, where we used that
\begin{equation*} 
  e^{\IF}= e^{\IF^0}\cdot(1+\sum_{i,j} t_{i,j}\IF^1_{i,j}) + o(t^2). 
\end{equation*}

{\it Proof of the theorem:} Instead of proving the master equation for the full descendant potential $\IF$, we first show that it suffices to prove 
\begin{equation*} 
   e^{\IF^0}\overleftarrow{\IH^{1,+}_{i,j}} - \overrightarrow{\IH^{1,-}_{i,j}}e^{\IF^0} = 
   \overrightarrow{\IH^{0,-}}(e^{\IF^0}\IF^1_{i,j}) - (e^{\IF^0}\IF^1_{i,j})\overleftarrow{\IH^{0,+}}. 
\end{equation*} 
Indeed, it is easy to see that the desired identity implies that
\begin{equation*} 
    F^{0,+}(\IH^{1,+}_{i,j}) - F^{0,-}(\IH^{1,-}_{i,j}) = 
    e^{+\IF^0}\overleftarrow{\IH^{1,+}_{i,j}}e^{-\IF^0} - e^{-\IF^0}\overrightarrow{\IH^{1,-}_{i,j}}e^{+\IF^0} 
\end{equation*}
is equal to  
\begin{equation*}
  e^{-\IF^0}\overrightarrow{\IH^{0,-}}(e^{+\IF^0}\IF^1_{i,j}) - (e^{+\IF^0}\IF^1_{i,j})\overleftarrow{\IH^{0,+}}e^{-\IF^0} 
  = D^{\IF^0}(\IF^1_{i,j}),
\end{equation*}
so that, after passing to homology, we have
\begin{equation*}  F^{0,+}_*(\IH^{1,+}_{i,j}) = F^{0,-}_*(\IH^{1,-}_{i,j}) \in H_*(\hbar^{-1}\DD^0,D^{\IF^0}) \end{equation*}
as desired. \\

On the other hand, the above identity directly follows from our definition of gravitational descendants of moduli spaces based on the definition of coherently connecting sections in tautological line bundles and the compactness theorem in [BEHWZ]. Indeed, in the same way as it is shown in [EGH] that the master equation for $\IF^0$ and $\IH^{0,\pm}$ follows from the fact that the codimension-one boundary of every moduli space $\CM_0$ is formed by products of moduli spaces $\CM_{1,0}\times\CM^+_{2,0}$ and $\CM^-_{1,0}\times\CM_{2,0}$, the desired identity relating $\IF^0$, $\IF^1_{i,j}$ and $\IH^{0,\pm}$, $\IH^{1,\pm}_{i,j}$ can be seen to follow from the fact that the codimension-one boundary of a descendant moduli space $\CM^j_1$ is given by products of the form $\CM^j_{1,1}\times\CM^+_{2,0}$, $\CM^{j,-}_{1,1}\times\CM_{2,0}$ and $\CM_{1,0}\times\CM^{j,+}_{2,1}$, $\CM^-_{1,0}\times\CM^j_{2,1}$: While the two summands involving $\IF^0$ and $\IH^{1,-}_{i,j}$, $\IH^{1,+}_{i,j}$ on the left-hand-side of the equation collect all boundary components of the form $\CM^{j,-}_{1,1}\times\CM_{2,0}$, $\CM_{1,0}\times\CM^{j,+}_{2,1}$, the two summands involving $\IF^1_{i,j}$ and $\IH^{0,-}$, $\IH^{0,+}$ on the right-hand-side of the equation collect all boundary components of the form $\CM^-_{1,0}\times\CM^j_{2,1}$, $\CM^j_{1,1}\times\CM^+_{2,0}$, respectively. Note that as for the master equation for $\IF^0$ and $\IH^{0,\pm}$ the appearance of $\IF^0$ in the exponential follows from the fact that there corresponding curves may appear with an arbitrary number of connected components, while the curves counted for in $\IH^{0,\pm}$, $\IH^{1,\pm}_{i,j}$, $\IF^1_{i,j}$ can only appear once due to index reasons or since there is just one additional marked point. \\

Finally, in order to see why we actually have $\IH^{1,-}_{i,j}=\IH^{1,+}_{i,j}$ on homology if we fixed $\lambda^-=\lambda^+=\lambda$, $\Ju^-=\Ju^+=\Ju$ and the abstract polyfold perturbations to have well-defined moduli spaces, observe that in this case $\IF^0$ just counts orbit cylinders, so that $F^{0,\pm}$ and hence $(F^{0,\pm})_*$ is the identity. $\qed$

\subsection{The circle bundle case}

In this subsection we briefly want to discuss the important case of circle bundles over closed symplectic manifolds, which links our constructions to gravitational descendants in Gromov-Witten theory, see also [R]. \\

For this recall that to any closed symplectic manifold $(M,\omega)$ with integral symplectic form $[\omega]\in H^2(M,\IZ)$ one can canonically assign a principal circle bundle $\pi: V\to M$ over $(M,\omega)$ by requiring that $c_1(V)=[\omega]$. Furthermore, it is easy to see that an $S^1$-connection form $\lambda$ with curvature $\omega$ on $\pi:V\to M$ is a contact form on the total space $V$, where the underlying contact structure agrees with the corresponding horizontal plane field $\xi=\ker\lambda$, while the Reeb vector field $R$ agrees with the infinitesimal generator of the $S^1$-action. Observe that a $\omega$-compatible almost complex structure $J$ on $M$ naturally equips $\IR\times V$ with a cylindrical almost complex structure by requiring that $\Ju$ maps the Reeb vector field to the $\IR$-direction and agrees with $J$ on the horizontal plane field $\xi$, which is naturally identified with $TM$. \\

Since every fibre of the circle bundle is hence a closed Reeb orbit for the contact form $\lambda$, it follows that the space of orbits is given by $M\times\IN$, where the second factor just refers to the multiplicity of the orbit. Hence, while every contact form in this class is not Morse as long as the symplectic manifold is not a point, it is still of Morse-Bott type. \\

Following [EGH] the Weyl algebra $\WW^0$ in this Morse-Bott case is now generated by sequences of graded formal variables $p_{\alpha,k}$, $q_{\alpha,k}$, $k\in\IN$ assigned to cohomology classes $\alpha$ forming a basis of $H^*(M,\IZ)$. For circle bundles in the Morse-Bott setup we now show that the general theorem from above leads the following stronger statement. Note that in the following theorem we do {\it not} assume that the sequences of coherent collections of sections are neccessarily $S^1$-invariant. \\
\\
{\bf Theorem 1.11:} {\it For circle bundles over symplectic manifolds, which are equipped with $S^1$-invariant contact forms, cylindrical almost complex structures (and abstract polyfold perturbations) as described above, the descendant Hamiltonians $\IH^1_{i,j}$ define a sequence of commuting operators on $\WW^0$, which is independent of the auxiliary data.} \\
\\
{\it Proof:} Observing that a map $\tilde{u}:(\Sigma,j)\to (\IR\times V,\Ju)$ from a punctured Riemann sphere to the cylindrical manifold $\IR\times V$, which is equipped with the canonical cylindrical almost complex structure $\Ju$ defined by the $\omega$-compatible almost complex structure $J$ on $M$, can be viewed as tuple $(h,u)$, where $u:(\Sigma,j)\to(M,J)$ is a $J$-holomorphic curve in $M$ and $h$ is a holomorphic section in $\IR\times u^*V\to\Sigma$, it is easy to see that every moduli spaces studied in SFT for the contact manifold $V$ carries a natural circle bundle structure after quotienting out the natural $\IR$-action. It follows that $D^0=0$, so that by our first theorem the $\IH^1_j$ already commute as elements in $\WW^0$. On the other hand, as long as the two different collections of auxiliary structures for $V$ are actually obtained as pull-backs of the corresponding auxiliary structures on $M$, it follows in the same way that the only rigid holomorphic curves in the resulting cobordisms are the orbit cylinders, so that the resulting automorphism is indeed the identity. $\qed$ \\   
 
For $S^1$ and $S^3$ Eliashberg already pointed out in his ICM 2006 talk, see [E], that the corresponding sequences $\Ih^1_j$ counting only genus zero curves lead to classical integrable systems, while the sequences of commuting operators $\IH^1_j$ provide deformation quantizations for these hierarchies. This is based on the surprising fact that the sequence $\Ih^1_j$ of Poisson-commuting functions actually agrees with integrable system for genus zero from Gromov-Witten theory obtained using the underlying Frobenius manifold structure. In particular, for $V=S^1$ it follows that that the resulting system of Poisson-commuting functions are precisely the commuting integrals of the dispersionless KdV hierarchy,
\begin{equation*} \Ih^1_j = \oint_{S^1} \frac{u^{j+2}(x)}{(j+2)!}\,dx,\;\;u(x)=\sum_{n\in\IN} p_n\;e^{+2\pi inx}+ q_n\;e^{-2\pi inx}, \end{equation*} 
while in the case of the Hopf fibration $V=S^3$ over $M=S^2$ one arrives at the Poisson-commuting integrals of the continuous limit 
of the Toda lattice. \\

In order to see why in genus zero the SFT of the circle bundle $V$ is so closely related to the Gromov-Witten theory of its symplectic base $M$, we recall from the proof of the theorem that every $\Ju$-holomorphic curve $\tilde{u}$ can be identified with a tuple $(h,u)$, where $u$ is a $J$-holomorphic curve in $M$ and $h$ is a holomorphic section in $\IR\times u^*V\to\Sigma$, whose poles and zeroes correspond to the positive and negative punctures with multiplicities. Since the zeroth Picard group of $S^2$ is trivial and hence every degree zero divisor is indeed a principal divisor, it follows that for every map $u$ the space of sections is isomorphic to $\IC$ and hence that the SFT moduli space of $\Ju$-holomorphic curves in $\IR\times V$ is indeed a circle bundle over the corresponding Gromov-Witten moduli space of $J$-holomorphic curves in $M$. \\

While this explains the close relation of SFT of circle bundles and Gromov-Witten theory in the genus zero case, the non-triviality of the Picard group for nonzero genus implies that the relation gets much more obscure when we allow for curves of arbitrary genus. Indeed, while in the case of $V=S^1$ the sequence $\IH^1_j$ defined by counting curves of arbitrary genus in $\IR\times V$ leads to the deformation quantization of the dispersionless KdV hierarchy, in particular, a quantum integrable system, counting curves of all genera in the underlying symplectic manifold, that is, the point, leads by Witten's conjecture to the classical integrable system given by the full KdV hierarchy as proven by Kontsevich . \\

At the end of this subsection we again want to emphasize that the above statement crucially relies on the fact that $V$ is equipped with a $S^1$-invariant contact form, cylindrical almost complex structure and abstract polyfold perturbations. Assuming for the moment that the sequences of coherent collections of sections are also chosen to be $S^1$-invariant, note that in this case the above invariance statement can directly be deduced from the independence of the descendant Gromov-Witten potential of the auxiliary data used to define it, which essentially relies on the fact that all moduli spaces have only boundary components of codimension greater or equal to two, so that absolute rather than relative virtual classes are defined. In particular, the gravitational descendants can be defined by integrating powers of the first Chern class over the absolute moduli cycle. On the other hand, recall that for the above theorem we did {\it not} require that the sequences of coherent collections of sections are neccessarily $S^1$-invariant. While our definition of coherent collections of sections seems to be very weak, our above theorem shows that the nice invariance property continues to hold even for a larger class of sections.

\section{Example: Symplectic field theory of closed geodesics}  

\subsection{Symplectic field theory of a single Reeb orbit}
We are now going to consider a concrete example, which actually formed the starting point for the formal discussion from above. \\

As above consider a closed contact manifold $V$ with chosen contact form $\lambda\in\Omega^1(V)$ and let $\Ju$ be a compatible cylindrical almost complex structure on $\IR\times V$. For any closed orbit $\gamma$ of the corresponding Reeb vector field $R$ on $V$ the 
orbit cylinder $\IR\times\gamma$ together with its branched covers are the basic examples of $\Ju$-holomorphic curves in $\IR\times V$. \\
 
In [F2] we prove that these {\it orbit curves} do not contribute to the algebraic invariants of symplectic field theory as long as they do not carry additional marked points. Our proof explicitly uses that the orbit curves (over a fixed orbit) are closed under taking boundaries and gluing, which follows from the fact that orbit curves are also trivial in the sense that they have trivial contact area and that this contact area is preserved under taking boundaries and gluing. In particular, it follows, see [F2], that every algebraic invariant of symplectic field theory has a natural analog defined by counting only orbit curves. Further specifying the underlying Reeb orbit let us hence introduce the {\it symplectic field theory of the Reeb orbit $\gamma$:} \\     

For this denote by $\WW^0_{\gamma}$ be the graded Weyl subalgebra of the Weyl algebra $\WW$, which is generated only by those $p$- and 
$q$-variables $p_n=p_{\gamma^n}$, $q_n=q_{\gamma^n}$ corresponding to Reeb orbits which are multiple covers of the fixed orbit $\gamma$ {\it and which are good in the sense of [BM]}. In the same way we further introduce the Poisson subalgebra $\PP^0_{\gamma}$ of $\PP^0$. It will become important that the natural identification of the formal variables $p_n$ and $q_n$ does {\it not} lead to an isomorphism of the graded algebras $\WW^0_{\gamma}$ and $\PP^0_{\gamma}$ with the corresponding graded algebras $\WW^0_{S^1}$ and $\PP^0_{S^1}$ for $\gamma=V=S^1$, not only since the gradings of $p_n$ and $q_n$ are different and hence even the commutation rules may change but also that variables $p_n$ and $q_n$ may not be there since they would correspond to bad orbits. \\

In the same way as we introduced the (rational) Hamiltonian $\IH^0$ and $\Ih^0$ as well as sequences of descendant Hamiltonians $\IH^1_j$ and $\Ih^1_j$ by counting general curves in the symplectization of a contact manifold, we can define distinguished elements $\IH^0_{\gamma}\in\hbar^{-1}\WW^0_{\gamma}$ and $\Ih^0_{\gamma}\in\PP^0_{\gamma}$, as well as sequences of descendant Hamiltonians $\IH^1_{\gamma,j}$ and $\Ih^1_{\gamma,j}$ by just counting branched covers of the orbit cylinder over $\gamma$ with signs 
(and weights), where the preservation of the contact area under splitting and gluing of curves proves that for every theorem from above we have a version for $\gamma$. \\

While for the general part described above we have already emphasized that the theorems are not yet theorems in the strict mathematical sense since the neccessary transversality theorems for the Cauchy-Riemann operator are part of the on-going polyfold project by Hofer and his collaborators and we further used the assumption that all occuring moduli spaces are manifolds with corners, for the rest of this paper we will {\it restrict to the rational case}, i.e., we will only be interested in the Poisson-commuting sequences $\Ih^1_{\gamma,j}$ on $H_*(\PP^0_{\gamma},d^0_{\gamma})$, but in return solve the occuring analytical problems in all detail. In particular, we have already proven in the paper [F2] that for (rational) orbit curves the transversality problem can indeed be solved using finite-dimensional obstruction bundles instead of infinite-dimensional polybundles. In order to see why this is even neccessary, observe that while in the case when $\gamma=V=S^1$ the Fredholm index equals the dimension of the moduli space, for general $\gamma\subset V$ the Fredholm index of a true branched cover is in general strictly smaller than the dimension of the moduli space of branched covers, so that transversality for the Cauchy-Riemann operator can in general not be satisfied. \\

So let us recall the main results about obstruction bundle transversality for orbit curves, where we refer to [F2] for all details. The first observation for orbit curves is that the cokernels of the linearized Cauchy-Riemann operators indeed fit together to give a smooth vector bundle $\overline{\Coker} \CR_{\Ju}$ over the compactified (nonregular) moduli spaces $\CM$ of orbit curves (of constant rank). It follows that every transveral section $\bar{\nu}$ of this cokernel bundle leads to a compact perturbation making the Cauchy-Riemann operator transversal to the zero section in the underlying polyfold setup. \\

In Gromov-Witten theory we would hence obtain the contribution of the regular perturbed moduli space by integrating the Euler class of the finite-dimensional obstruction bundle over the compactified moduli space. On the other hand, passing from Gromov-Witten theory back to symplectic field theory again, we see that we just arrive at the same problem we had to face with when we wanted to define gravitational descendants in symplectic field theory. Indeed, as for the tautological line bundles, the presence of codimension-one boundary of the (nonregular) moduli spaces of branched covers implies that Euler numbers for sections in the cokernel bundles are not defined in general, since the count of zeroes depends on the compact perturbations chosen for the moduli spaces in the boundary. \\

Instead of looking at a single moduli space, we hence again have to consider all moduli spaces at once. Replacing the tautological line bundle $\LL$ by the cokernel bundle $\overline{\Coker}\CR_{\Ju}$ and considering the nonregular moduli space of branched covers instead of the regular moduli space itself, we hence now define {\it coherent} collections of sections in the obstruction bundles $\overline{\Coker}\CR_J$ over all moduli spaces $\CM$ as follows. \\

Following the compactness statement in [BEHWZ] for the contact manifold $S^1$ the codimension-one boundary of every moduli space of branched covers $\CM$ again consists of curves with two levels (in the sense of [BEHWZ]), whose moduli spaces can be represented as products $\CM_1\times\CM_2$ of moduli spaces of strictly lower dimension, where the first index again refers to the level. On the other hand, it follows from the linear gluing result in [F2] that over the boundary component $\CM_1\times\CM_2$ the cokernel bundle $\overline{\Coker}\CR_{\Ju}$ is given by 
\begin{equation*} \overline{\Coker}\CR_{\Ju}|_{\CM_1\times\CM_2} = \pi_1^*\overline{\Coker}^1\CR_{\Ju}\oplus\pi_2^*\overline{\Coker}^2\CR_{\Ju}, \end{equation*}
where $\overline{\Coker}^1\CR_{\Ju}$, $\overline{\Coker}^2\CR_{\Ju}$ denotes the cokernel bundle over the moduli space $\CM_1$, $\CM_2$ and $\pi_1$, $\pi_2$ is the projection onto the first or second factor, respectively. \\

Assuming that we have chosen sections $\bar{\nu}$ in the cokernel bundles $\overline{\Coker}\CR_{\Ju}$ over all moduli spaces $\CM$ of branched covers, we again call this collection of sections $(\bar{\nu})$ {\it coherent} if over every codimension-one boundary component $\CM_1\times\CM_2$ of a moduli space $\CM$ the corresponding section $\bar{\nu}$ agrees with the pull-back $\pi_1^*\bar{\nu}_1\oplus\pi_2^*\bar{\nu}_2$ of the chosen sections $\bar{\nu}_1$, $\bar{\nu}_2$ in the cokernel bundles $\overline{\Coker}^1\CR_{\Ju}$ over $\CM_1$, $\overline{\Coker}^2\CR_{\Ju}$ over $\CM_2$, respectively.  \\

Since in the end we will again be interested in the zero sets of these sections, we will again assume that all occuring sections are transversal to the zero section.
As before it is not hard to see that one can always find such coherent collections of (transversal) sections in the cokernel bundles by using induction on the dimension of the underlying nonregular moduli space of branched covers. Note that the latter is {\it not} equal to the Fredholm index. \\

In [F2] we prove the following result about orbit curves with {\it no} additional marked points. \\
\\
{\bf Theorem ([F2]):} {\it For the cokernel bundle $\overline{\Coker}\CR_J$ over the compactification $\CM$ of every moduli space 
of branched covers over an orbit cylinder with $\dim\IM-\rank\Coker\CR_J=0$ the following holds:}
\begin{itemize} 
\item {\it For every pair $\bar{\nu}^0$, $\bar{\nu}^1$ of coherent and transversal sections in $\overline{\Coker}\CR_J$ the algebraic count of zeroes 
      of $\bar{\nu}^0$ and $\bar{\nu}^1$ are finite and agree, so that we can define an Euler number 
      $\chi(\overline{\Coker}\CR_J)$ for coherent sections in $\overline{\Coker}\CR_J$ by} 
      \begin{equation*} \chi(\overline{\Coker}\CR_J) \,:=\, \sharp (\bar{\nu}^0)^{-1}(0) \,=\, \sharp (\bar{\nu}^1)^{-1}(0). 
      \end{equation*} 
\item {\it This Euler number is $\chi(\overline{\Coker}\CR_J) = 0$.} 
\end{itemize}
$ $\\
This theorem in turn has the following consequence.\\
\\
{\bf Corollary 2.1:} { \it For every closed Reeb orbit $\gamma$ the Hamiltonian $\Ih^0_{\gamma}$ vanishes independently of the chosen coherent collection of sections $(\bar{\nu})$ in the cokernel bundles over all moduli spaces of branched covers,}
\begin{equation*} \Ih^0 = \Ih^{0,\bar{\nu}}= 0. \end{equation*}
{\it In particular, the sequences of descendant Hamiltonians $\Ih^1_{\gamma,j}$ already Poisson-commute as elements in $\PP^0_{\gamma}$.}\\

Note that the latter statement is obvious in the case $\gamma=V=S^1$. While it directly follows from index reasons that $\Ih^1_{S^1,j}=0$ when the string of differential forms just consists of the zero-form $1$ on $S^1$, it is shown in [E] using the results from Okounkov and Pandharipande in [OP] that for the one-form $dt$ on $S^1$ the system of Poisson commuting functions on $\PP^0_{S^1}$ is given by
\begin{equation*} \Ih^1_{S^1,j} = \oint_{S^1} \frac{u^{j+2}(x)}{(j+2)!}\,dx,\;\;u(x)=\sum_{n\in\IN} p_n\;e^{+2\pi inx}+ q_n\;e^{-2\pi inx}, \end{equation*} 
i.e., hence agrees with the {\it dispersionless KdV (or Burger) integrable hierarchy.} \\
 
Going back from $\gamma=V=S^1$ to the case of orbit curves over general Reeb orbits $\gamma$, observe that, since for the orbit curves the evaluation map to $V$ factors through the inclusion map $\gamma\subset V$, it follows that it again only makes sense to consider zero- or one-forms, where we can assume without loss of generality that the zero-form agrees with $1\in\Omega^0(V)$ and that the integral of the one-form $\theta\in\Omega^1(V)$ over the Reeb orbit is one,
\begin{equation*} \int_{\gamma}\theta = 1. \end{equation*}
For the case with no gravitational descendants, note that it follows from index reasons that the only curves to be considered are orbit cylinders with one marked point, since introducing an additional marked point adds two or one to the Fredholm index. Since orbit cylinders are always regular and their contribution hence just equals the integral of the form $\theta$ over the closed orbit $\gamma$, we hence get just like in the case of $\gamma=V=S^1$ that the zeroth descendant Hamiltonian $\Ih^1_{\gamma,0}$ vanishes if $\deg\theta = 0$ and 
\begin{equation*} \Ih^1_{\gamma,0} = \oint_{S^1} \frac{u^2(x)}{2!}\,dx = \sum p_n q_n \end{equation*}
if $\deg\theta = 1$ with the normalization from above. For the sum note that we only assigned formal variables $p_n$, $q_n$ to Reeb orbits which are good in the sense of [BM]. \\

While the Hamiltonians $\Ih^1_{\gamma,0}$ hence agree with the Hamiltonian $\Ih^1_{S^1,0}$ for $\gamma=V=S^1$ up to the problem of bad orbits, since no obstruction bundles have to be considered, it is easy to see that the argument breaks down when gravitational descendants are introduced, since the underlying orbit curve then has non-zero Fredholm index $1+2(j-1)+\deg \theta$ and hence need not be an orbit cylinder anymore. While for the case of a one-form we can hence expect to find new integrals for the nontrivial Hamiltonian $\Ih^1_{S^1,0}=\Ih^1_{\gamma,j}$, we first show that in the case of a zero-form not only the zeroth Hamiltonian but even the whole sequence of descendant Hamiltonians $\Ih^1_{\gamma,j}$ is trivial.  \\
\\  
{\bf Theorem 2.2:} {\it Let $\gamma$ be a Reeb orbit in any contact manifold $V$ and assume that the string of differential forms on $V$ just consists of the zero-form $1\in\Omega^0(V)$. Then the sequence of Poisson-commuting functions $\Ih^1_{\gamma,j}$ on $\PP^0_{\gamma}$ is trivial,}
\begin{equation*} \Ih^1_{\gamma,j} = 0,\;j\in\IN \end{equation*}
{\it just like in the case of $\gamma=V=S^1$.}\\
\\
{\it Proof:} Since the proof of this theorem follows from completely the same arguments as the proof of our theorem in [F2] about Euler numbers of coherent sections in obstruction bundles from above, we shortly give the main idea for the proof in [F2] about orbit curves without additional marked points and then discuss its generalization to orbit curves with zero-forms and gravitational descendants. \\
   
After proving that we can work with finite-dimensional obstruction bundles instead of infinite-dimensional polybundles, recall that the main problem lies in the presence of codimension-one boundary of the (nonregular) moduli space, so that Euler numbers of Fredholm problems are not defined in general, since the count of zeroes in general depends on the compact perturbations chosen for the moduli spaces in the boundary. In [F2] we prove the existence of the Euler number for moduli spaces of orbit curves without additional marked points by induction on the number of punctures. For the induction step we do not only use that there exist Euler numbers for the moduli spaces in the boundary, but it is further important that all these Euler numbers are in fact trivial. The vanishing of the Euler number in turn is deduced from the different parities of the Fredholm index of the Cauchy-Riemann operator and the actual dimension of the moduli space of branched covers following the idea for the vanishing of the Euler characteristic for odd-dimensional manifolds. \\

For the generalization to the case of additional marked points and gravitational descendants, it is clear that it still suffices to work with finite-dimensional obstruction bundles. On the other hand, recall that the only further ingredient to our proof in [F2] was that the Fredholm index and the dimension of the moduli spaces always have different parity. Hence it follows that the proof in [F2] also works for the case when $\theta$ is a zero-form as the actual dimension of the moduli spaces is still even, while it breaks down in the case when $\theta$ is a one-form. $\qed$ \\

Observe that for one-forms it is indeed no longer clear that the every Euler number has to be zero, as we for $\gamma=V=S^1$ and $\theta=dt$ we get nontrivial contributions from true branched covers. While at first glance the major problem seems to be the truely complicated computation of the Euler number (see [HT1], [HT2] for related results), we further have the problem that Euler numbers need no longer exist for all Fredholm problems. {\it For the rest of this paper we will hence only be interested in the case where the chosen differential form has degree one,} $\deg\theta =1$.\\

While for $\gamma=V=S^1$ we actually get a unique sequence of Poisson-commuting functions, observe that for general fixed Reeb orbits $\gamma$ in contact manifolds $V$ the descendant Hamiltonians $\Ih^1_{\gamma,j} = \Ih^{1,\bar{\nu}}_{\gamma,j}$ may indeed depend on the chosen collection of sections in the cokernel bundles $\overline{\Coker}\CR_{\Ju}$. Hence the invariance statement is no longer trivial, but implies that for different choices of coherent abstract perturbations $\bar{\nu}^{\pm}$ for the moduli spaces the resulting system of commuting elements $\Ih^{1,-}_{\gamma,j}$, $j=0,1,2,..$ and $\Ih^{1,+}_{\gamma,j}$, $j=0,1,2,..$ on $\PP_{\gamma}^0$ are just isomorphic, i.e., there exists an {\it auto}morphism of the Poisson algebra $\PP_{\gamma}^0$ which identifies $\Ih^{1,-}_{\gamma,j}\in \PP_{\gamma}^0$ with $\Ih^{1,+}_{\gamma,j}\in \PP_{\gamma}^0$ for all $j\in\IN$. \\

The above discussion hence shows that the computation of the symplectic field theory of a closed Reeb orbit gets much more difficult when gravitational descendants are considered. In what follows we want to determine it in the special case where the contact manifold is the unit cotangent bundle $S^*Q$ of a ($m$-dimensional) Riemannian manifold $Q$, so that every closed Reeb orbit $\gamma$ on $V=S^*Q$ corresponds to a closed geodesic $\bar{\gamma}$ on $Q$. \\

Before we can state the theorem we first want to expand the descendant Hamiltonians $\Ih^1_{S^1,j}$ in terms of the $p_n$- and $q_n$-variables, where set $p_n=q_{-n}$. Abbreviating $u_n(x) = q_n e^{inx}$ for every nonzero integer $n$ it follows from $u=\sum_n u_n$ that 
\begin{equation*} 
 \Ih^1_{S^1,j} = \oint_{S^1} \frac{u^{j+2}(x)}{(j+2)!}\,dx
               = \oint_{S^1} \sum \frac{u_{n_1}(x)\cdot ... \cdot u_{n_{j+2}}(x)}{(j+2)!}\, dx
\end{equation*}
On the other hand, note that the integration around the circle corresponds to selecting only those sequences of multiplicities $(n_1,...,n_{j+2})$, whose sum is equal to zero, so that 
\begin{eqnarray*} 
\Ih^1_{S^1,j} \,=\, \sum_{n_1+...+n_{j+2} = 0} \frac{q_{n_1}\cdot ... \cdot q_{n_{j+2}}}{(j+2)!}.
\end{eqnarray*} 
Apart from the sequence of Poisson-commuting functions for the circle, the grading of the functions given by the grading of $p_n$- and $q_n$-variables will play a central role for the upcoming theorem. For this observe that it follows from the grading conventions in symplectic field theory that the grading of the full Hamiltonian $\IH^0$ is $-1$, so that by $\IH^0 = \sum_g \hbar^{g-1} \IH^0_g$ the grading for the rational Hamiltonian $\Ih^0=\IH^0_0$ is given by $|\Ih^0| = |\IH^0|+|\hbar| = -1+2(m-2)$. Since this grading has to agree with the grading of $t_j\Ih^1_j$ with $|t_j|=2(1-j)-\deg\theta = 1-2j$, it follows that for every Reeb orbit $\gamma\subset V$ we have
\begin{equation*} |\Ih^1_{\gamma,j}| = -1+2(m-2)-1+2j = 2(m+j-3). \end{equation*}

We already mentioned that the natural identification of the formal variables $p_n$ and $q_n$ does {\it not} lead to an isomorphism of the graded algebras $\WW^0_{\gamma}$ and $\PP^0_{\gamma}$ with the corresponding graded algebras $\WW^0_{S^1}$ and $\PP^0_{S^1}$ for $\gamma=V=S^1$, not only since the gradings of $p_n$ and $q_n$ are different and hence even the commutation rules may change but even that variables $p_n$ and $q_n$ may not be there since they would correspond to bad orbits. While for the grading of $\gamma=V=S^1$ given by $|p_n|=|q_n|=-2$ in the descendant Hamiltonians $\Ih^1_{S^1,j}$ every summand indeed has the same degree $2(m+j-3)$, passing over to a general Reeb orbit $\gamma$ with the new grading given by $|p_n|=m-3-\CZ(\gamma^n)$, $|q_n|=m-3+\CZ(\gamma^n)$ the descendant Hamiltonian $\Ih^1_{S^1,j}$ is no longer of pure degree, i.e., different summands of the same descendant Hamiltonian usually have different degree. While the Poisson-commuting sequence for the circle seems not to be related to the sequence of descendant Hamiltonians for general Reeb orbits $\gamma$, we prove the following result in the case when the Reeb orbit corresponds to a closed geodesic. \\
\\
{\bf Theorem 2.3:} {\it Assume that the contact manifold is the unit cotangent bundle $V=S^*Q$ of a Riemannian manifold $Q$, so that the closed Reeb orbit $\gamma$ corresponds to a closed geodesic $\bar{\gamma}$ on $Q$, and that the string of differential forms just consists of a single one-form which integrates to one around the orbit. Then the resulting system of Poisson-commuting functions $\Ih^1_{\gamma,j}$, $j\in\IN$ on $\PP^0_{\gamma}$ is isomorphic to the system of Poisson-commuting functions $\Ig^1_{\bar{\gamma},j}$, $j\in\IN$ on $\PP^0_{\bar{\gamma}}=\PP^0_{\gamma}$, where for every $j\in\IN$ the descendant Hamiltonian $\Ig^1_{\bar{\gamma},j}$ is given by} 
\begin{equation*} 
 \Ig^1_{\bar{\gamma},j} \;=\; \sum \epsilon(\vec{n})\frac{q_{n_1}\cdot ... \cdot q_{n_{j+2}}}{(j+2)!} 
\end{equation*}
{\it where the sum runs over all ordered monomials $q_{n_1}\cdot ... \cdot q_{n_{j+2}}$ with $n_1+...+n_{j+2} = 0$ \textbf{and which are of degree $2(m+j-3)$}. Further $\epsilon(\vec{n})\in\{-1,0,+1\}$ is fixed by a choice of coherent orientations in symplectic field theory and is zero if and only if one of the orbits $\gamma^{n_1},...,\gamma^{n_{j+2}}$ is bad.} \\

We have the following immediate corollary, which immediately follows from the behavior of the Conley-Zehnder index for multiple covers.  \\
\\
{\bf Corollary 2.4:} {\it Assume that the closed geodesic $\bar{\gamma}$ represents a hyperbolic Reeb orbit in the unit cotangent bundle and $\dim Q>1$. Then $\Ig^1_{\bar{\gamma},j}=0$ and hence $\Ih^1_{\gamma,j}=0$ for all $j>0$.} \\

Indeed, since for hyperbolic Reeb orbits the Conley-Zehnder index $\CZ(\gamma^n)$ of $\gamma^n$ is given by $\CZ(\gamma^n)=n\cdot\CZ(\gamma)$, an easy computation shows that there are {\it no} products of the above form of the desired degree. On the other hand, note that without the degree condition we would just get back the sequence of descendant Hamiltonians for the circle. Forgetting about orientation issues, in simple words we can hence say that the sequence $\Ig^1_{\bar{\gamma},j}$ is obtained from the sequence for $\bar{\gamma}=Q=S^1$ by removing all summands with the wrong, that is, not maximal degree, where the latter can explicitely be computed using the formulas in [Lo] but also follows from our proof. \\

The proof relies on the observation that for orbit curves the gravitational descendants indeed have a geometric meaning in terms of branching conditions, which is a slight generalization of the result for the circle shown by Okounkov and Pandharipande in [OP]. Applying (and generalizing) the ideas of Cieliebak and Latschev in [CL] for relating the symplectic field theory of $V=S^*Q$ to the string topology of the underlying Riemannian manifold $Q$, we then study branched covers of the corresponding trivial half-cylinders in the cotangent bundle connecting the Reeb orbit $\gamma$ with the underlying geodesic $\bar{\gamma}$ to prove that the sequence of Poisson-commuting functions $\Ih^1_{\gamma,j}$ is isomorphic to a sequence of Poisson-commuting functions $\Ig^1_{\bar{\gamma},j}$. {\it While the descendant Hamiltonians $\Ih^1_{\gamma,j}$ on the SFT side are defined using very complicated obstruction bundles over (nonregular) moduli spaces of arbitary large dimension, the key observation is that for the descendant Hamiltonians $\Ig^1_{\bar{\gamma},j}$ on the string side we indeed only have to study obstruction bundles over discrete sets, which clearly disappear if the Fredholm index is right.} With this we get that the Poisson-commuting sequences for the closed geodesics can be computed from the sequences for the circle {\it and} the Morse indices of the geodesic and its iterates as stated in the theorem. \\

\subsection{Gravitational descendants = branching conditions}
Recall that by the above theorem from the last subsection we only have to consider the case where $\theta$ is a one-form on $V$, where we still assume without loss of generality that the integral of $\theta$ over $\gamma$ is one. It follows that integrating the pullback of $\theta$ under the evaluation map over the moduli space of orbit curves with one additional marked point and dividing out the natural $\IR$-action on the target $\RS\cong\IR\times\gamma$ is equivalent to restricting to orbit curves where the additional marked point is mapped to a special point on $\RS$. In other words, in what follows we will view $h^1_{\gamma,j}$ no longer as part of the Hamiltonian for $\gamma$ but as part of the potential for the cylinder over $\gamma$ equipped with a non-translation-invariant two-form. In order to save notation, $\CM_1=\CM_1(\Gamma^+,\Gamma^-)$ will from now on denote the corresponding moduli space. On the other hand, after introducing coherent collections $(\bar{\nu})$ of obstruction bundle sections, it is easy to see that the tautological line bundle $\LL^{\bar{\nu}}$ over $\CM_1^{\bar{\nu}}$ is just the restriction of the tautological line bundle $\LL$ over $\CM_1$ to $\CM_1^{\bar{\nu}}=\bar{\nu}^{-1}(0)\subset\CM_1$. \\

For the orbit curves we now want to give a geometric interpretation of gravitational descendants in terms of branching conditions over the special point on $\RS$. Before we state the corresponding theorem and give a rigorous proof using the stretching-of-the-neck procedure from SFT, we first informally describe a naive direct approach based on our definition of gravitational descendants from above, which should illuminate the underlying geometric ideas. \\

Recall that if $(h,z)$ is an element in the non-compactified moduli space $\IM^{\nu}_1\subset\CM^{\bar{\nu}}_1$ the fibre of the canonical line bundle $\LL$ over $(h,z)$ is given by $\LL_{(h,z)}=T_z^*S$. Identifying the tangent space to the cylinder at the special point with $\IC$ it follows that $s(h,z)=\frac{\del h}{\del z}(z)\in T_z^*S$ is a section in the restriction of $\LL$ to $\IM^{\nu}_1$. Since $s$ is a transversal section in the tautological line bundle over $\IM_1^{\nu}$ if and only if it extends to a section over $\IM_1$ such that $s\oplus\nu$ is transversal to the zero section in $\LL\oplus\Coker\CR_{\Ju}$ over $\IM_1$, we may assume after possibly perturbing $\nu$ that $s$ is indeed transversal. On the other hand, since $\frac{\del h}{\del z}(z)=0$ is equivalent to saying that $z\in S$ is a branch point of the holomorphic map $h:S\to \CP$, it follows that $\IM_1^1:=s^{-1}(0)\subset\IM_1$ indeed agrees with the space of all orbit curves $(h,z)$ with one additional marked point, where $z$ is a branch point of $h$. \\

Further moving on to the case $j=2$ observe that a natural candidate for a generic section $s_2$ in the restriction of the product line bundle $\LL^{\otimes 2}$ to $\IM^1_1\subset\IM_1$ is given by $s_2(h,z)=\frac{\del^2 h}{\del z^2}(z)\in (T_z^*S)^{\otimes 2}$, for which $\IM_1^2=s_2^{-1}(0)\subset\IM_1^1$ agrees with the space of holomorphic maps where $z\in S$ is now a branch point of order at least two. For general $j$ we can hence proceed by induction and define the section $s_j$ in $\LL^{\otimes j}$ over $\IM^{j-1}_1:=s_{j-1}^{-1}(0)\subset\IM_1$ by $s_j(h,z)=\frac{\del^j h}{\del z^j}(z)$, so that $\IM^j_1$ agrees with the space of holomorphic maps where $z\in S$ is a branch point of order at least $j$. \\

If the chosen sections $s_1,...,s_j$ over the non-compactified moduli spaces would extend in the same way to a coherent collection of sections in the tautological line bundles over the compactified moduli spaces $\CM_1^{\bar{\nu}}$, the above would show that in the case of orbit curves considering the $j$.th descendant moduli space is equivalent after passing to homology to requiring that the underlying additional marked point is a branch point of order $j$. In [OP] it was however shown that already for the case of the circle $\gamma=V=S^1$ the latter assumption is not entirely true, but that one instead additionally obtains corrections from the boundary $\CM_1-\IM_1$. \\

To this end, we define a {\it branching condition} to be a tuple of natural numbers $\mu=(\mu_1,...,\mu_{\ell(\mu)})$ of length $\ell(\mu)$ and total branching order $|\mu|=\mu_1+...+\mu_{\ell(\mu)}$. Then the moduli space $\CM^{\mu}=\CM^{\mu}(\Gamma^+,\Gamma^-)$ consists of orbit curves with $\ell(\mu)$ connected components, where every connected component carries one additional marked point $z_i$, which is mapped to the special point on $\IR\times\gamma$ and is a branch point of order $\mu_i-1$ for $i=1,...,\ell(\mu)$. For every branching condition $\mu=(\mu_1,...,\mu_{\ell(\mu)})$ we then define new Hamiltonians $\Ih^1_{\gamma,\mu}=\Ih^{1,\bar{\nu}}_{\gamma,\mu}$ by setting 
\begin{equation*} \Ih^1_{\gamma,\mu} = \sum_{\Gamma^+,\Gamma^-} \#\CM_1^{\mu}(\Gamma^+,\Gamma^-)\;p^{\Gamma^+}q^{\Gamma^-}. \end{equation*} 

With the following theorem we will prove that the abstract descendants-branching correspondence from [OP] holds for every closed Reeb orbit $\gamma\subset V$. 
For every $j\in\IN$ and every branching condition $\mu$ we let $\rho^0_{j,\mu}$ be the number given by integrating the $j$.th power of the first Chern class of the tautological line bundle over the moduli space of connected rational curves over $\CP$ with one marked point mapped to $0$ and $\ell(\mu)$ additional marked points $z_i$ mapped to $\infty$ which are branch points of order $\mu_i-1$, $i=1,...,\ell(\mu)$.  \\ 
\\
{\bf Lemma 2.5:} {\it Each of the descendant Hamiltonians $\Ih^1_{\gamma,j}$ can be written as a sum,} 
\begin{equation*} 
 \Ih^1_{\gamma,j} \;=\; \frac{1}{j!}\;\cdot\;\Ih^1_{\gamma,(j+1)} \;+\;\sum_{|\mu|<j} \rho^0_{j,\mu}\;\cdot\; \Ih^1_{\gamma,\mu}, 
\end{equation*} 
{\it where $\Ih^1_{\gamma,\mu}\in\PP^0_{\gamma}$ counts branched covers of the orbit cylinder with $\ell(\mu)$ connected components, where each component carries one additional marked point $z_i$, which is mapped to the special point on $\IR\times\gamma$ and is a branch point of order $\mu_i-1$ for $i=1,...,\ell(\mu)$.} \\

Note that the statement of the lemma can be rephrased by saying that the integration of the $j$.th power of the first Chern class corresponds to a weighted sum of branching conditions,
\begin{equation*} c_1(\LL)^j \;\doteq\; \frac{1}{j!}\;\cdot\;(j+1) \;+\; \sum_{|\mu|<j} \rho^0_{j,\mu} \;\cdot\; \mu. \end{equation*} 
which is the rational version of the abstract descendants-branching correspondence from [OP] for the circle, where the coefficient $\rho^0_{\mu,j}$ is nonzero and agrees with the coefficient $\rho_{j,\mu}$ from [OP] only if the genus $g$ determined by the Fredholm index,
\begin{equation*} j+1 = 2g-1+|\mu|+\ell(\mu) \end{equation*}
is zero. \\
\\
{\it Proof:} Recall that the result in [OP] for the circle relies on the degeneration formula from relative Gromov-Witten theory, where the target sphere with three special points $x^-=0$, $x^+=\infty$ and $x$ degenerates in such a way that the original sphere only carries the two special points $x^-=0$, $x^+=\infty$ while the third special point sits on a  second sphere connected to the original one by a node. Viewing a sphere with two special points as a cylinder, it is clear that a corresponding statement can be proven for a Reeb orbit $\gamma$ in a general contact manifold if the standard cylinder is replaced by the orbit cylinder $\IR\times\gamma$ in the symplectization of the contact manifold which degenerates to an orbit cylinder with a ghost bubble attached. \\

Since the degeneration formula from relative Gromov-Witten theory is no longer applicable, we will have to use the neck-stretching process from symplectic field theory, which however agrees with the degeneration process from relative Gromov-Witten theory in the case of the circle. For this observe that performing a neck-stretching at a small circle around the special point on the standard cylinder we obtain a pair-of-pants together with a complex plane carrying the special point, which can be identified with spheres with three or two special points, respectively. Replacing the circle by a Reeb orbit $\gamma$ in a general contact manifold the neck-stretching yields besides a complex plane with a special point a pair-of-pants with a positive and a negative cylindrical ends over $\gamma$ together with a cylindrical end over the circle. \\

Note that in order to include infinitesimal deformations needed for the obstruction bundles, we identify the orbit cylinder $\IR\times\gamma$ (together with an infinitesimal tubular neighborhood) with (an infinitesimal neighborhood of the zero section in) its normal bundle over $\RS$ with fibre given by the contact distribution $\xi$ and twist around the puncture given by the linearized Reeb flow along $\gamma$. Then the (infinitesimal) neck-stretching is performed along the (infinitesimal) hypersurface given by the restriction of the normal bundle to the small circle in $\RS$. \\

Before we make the proof rigorous by studying coherent collections of sections in the cokernel bundles and the tautological line bundles over the moduli space of branched covers for the circle, observe that theorem 2.5.5 in [EGH] concerning composition of cobordisms suggests that $\Ih^1_{\gamma,j}$, viewed as a potential on $\PP^0_{\gamma}$, is homotopic, and by $\Ih^0_{\gamma}=0$ hence agrees with a potential, which can directly be computed from the potential for the complex plane counting rational curves with one additional marked point mapped to the special point and the potential for the pair-of-pants with its cylindrical ends over $\gamma$ and the circle counting rational curves with no additional marked points. Indeed it follows from the compactness statement in [BEHWZ] that under the neck-stretching procedure every branched cover of the orbit cylinder with one additional marked point mapped to the special point splits into a branched cover of the complex plane with one additional marked point mapped to the special point and a branched cover of the pair-of-pants with no additional marked points. \\

While from $S^1$-symmetry reasons the potential for the complex plane with one special point can only count connected curves, note that under the splitting process the connected curve may split into branched covers of the pair-of-pants with more than one connected component. On the other hand, since the glued curve has genus zero, it follows that the branched cover of the complex plane and any connected component of the branched cover of the pair-of-pants cannot be glued at more than one cylindrical end, so that the number of connected components of the branched cover of the pair-of-pants 
agrees with the number of cylindrical ends of the branched cover of the complex plane. \\

Note that a collection $\Gamma$ of closed Reeb orbits in the contact manifold $S^1$ is naturally identified with a tuple $\mu=(\mu_1,...,\mu_{\ell(\mu)})$ of multiplicities and a branched cover is asymptotically cylindrical over the $\mu_i$.th iterate of the circle near the puncture $z_i$ precisely if $z_i$ is a branch point of order $\mu_i-1$. With this it follows that $\Ih^1_{\gamma,j}$ can be computed as desired by summing over all branching conditions $\mu=(\mu_1,...,\mu_{\ell(\mu)})$, where for each $\mu$ the summand is given by the product of $\rho^0_{j,\mu}$, obtained by integrating the $j$.th power of the first Chern class of the tautological line bundle over the moduli space of branched covers of $\CP$ with one marked point mapped to the special point and $\ell(\mu)$ additional marked points $z_i$ mapped to $\infty$ which are branch points of order $\mu_i-1$, with the branching Hamiltonians $\Ih^1_{\gamma,\mu}\in\PP^0_{\gamma}$, counting branched covers of the orbit cylinder with $\ell(\mu)$ connected components, where each component carries one additional marked point $z_i$, which is mapped to the special point on $\IR\times\gamma$ and is again a branch point of order $\mu_i-1$ for $i=1,...,\ell(\mu)$. \\

In order to make the proof rigorous it remains to understand the above statement on the level of coherent collections of sections in the cokernel bundles and the tautological line bundles over the moduli spaces of branched covers for the circle. For this observe that for every chosen collections of Reeb orbits $\Gamma^+,\Gamma^-$ the neck-stretching procedure at a small circle around the special point on the standard cylinder leads to a compactified moduli space $\widetilde{\IM}_1=\widetilde{\IM}_1(\Gamma^+,\Gamma^-)$. It is shown in [BEHWZ] that this compactified moduli space has the desired codimension-one boundary components $\CM_1=\CM_1(\Gamma^+,\Gamma^-)$ counting branched covers of the original orbit cylinder with one special point and $\CM_{1,1}\times\CM_{2,0}$ with $\CM_{1,1}=\CM_1(\Gamma)$ and $\CM_{2,0}=\CM_0(\Gamma^+,\Gamma,\Gamma^-)$ counting branched covers of the complex plane with one additional marked point and of the pair-of-pants with possibly more than one connected component, respectively. On the other hand, in contrast to the degeneration process from relative Gromov-Witten theory, it follows from [BEHWZ] that one also has to consider codimension-one boundary strata of the form $\widetilde{\IM}_{1,1}\times\CM_{2,0}$ with $\widetilde{\IM}_{1,1}=\widetilde{\IM}_1(\Gamma^+_1,\Gamma^-_1)$, $\CM_{2,0}=\CM_0(\Gamma^+_2,\Gamma^-_2)/\IR$ and $\CM_{1,0}\times \tilde{\IM}_{2,1}$ with $\CM_{1,0}=\CM_0(\Gamma^+_1,\Gamma^-_1)/\IR$, $\widetilde{\IM}_{2,1}=\widetilde{\IM}_1(\Gamma^+_2,\Gamma^-_2)$, which correspond to a splitting of a curve into two levels during the stretching process and which are irrelevant for the case of the circle due to the $S^1$-symmetry on $\CM_{1,0}$ and $\CM_{2,0}$, respectively. \\

First, since the coherent collections of sections in the cokernel bundles over the moduli spaces of branched covers by definition are {\it not} affected by the position of the additional marked point, it follows that one can use the same obstruction bundle perturbations $\bar{\nu}=\bar{\nu}(\Gamma^+,\Gamma^-)$ throughout the stretching process. In particular, it follows that the regular moduli space $\widetilde{\IM}_1^{\bar{\nu}}=\bar{\nu}^{-1}(0)\subset\widetilde{\IM}_1$ has codimension-one boundary components $\CM_1^{\bar{\nu}}$ and $\CM_{1,1}\times\CM_{2,0}^{\bar{\nu}}$ as well as $\widetilde{\IM}_{1,1}^{\bar{\nu}_1}\times\CM_{2,0}^{\bar{\nu}_2}$ and $\CM_{1,0}^{\bar{\nu}_1}\times \widetilde{\IM}_{2,1}^{\bar{\nu}_2}$, where $\bar{\nu}_1$, $\bar{\nu}_2$ are sections in the cokernel bundles $\overline{\Coker}^1\CR_{\Ju}$, $\overline{\Coker}^2\CR_{\Ju}$ over $\CM_{1,0}=\CM_0(\Gamma^+_1,\Gamma^-_1)$, $\CM_{2,0}=\CM_0(\Gamma^+_2,\Gamma^-_2)$ and which are determined by $\bar{\nu}$ by the coherency condition. \\

On the other hand, concerning the coherent collections of sections in the tautological line bundles, it can be shown as above that the tautological line bundle $\widetilde{\LL}=\widetilde{\LL}^{\bar{\nu}}$ over $\widetilde{\IM}_1^{\bar{\nu}}$ agrees with the tautological line bundle $\LL$ over $\CM_1^{\bar{\nu}}$, with the pullback $\pi_1^*\LL_1$ over $\CM_{1,1}\times\CM_{2,0}^{\bar{\nu}}$ and with the pullbacks $\pi_1^*\widetilde{\LL}_1$ over $\widetilde{\IM}_{1,1}^{\bar{\nu}_1}\times\CM_{2,0}^{\bar{\nu}_2}$ and $\pi_2^*\widetilde{\LL}_2$ over $\CM_{1,0}^{\bar{\nu}_1}\times \widetilde{\IM}_{2,1}^{\bar{\nu}_2}$, respectively. Assuming that we have chosen coherent collections of sections $(s)$ in the tautological line bundles $\LL$ over all moduli spaces $\CM_1^{\bar{\nu}}= \CM_1^{\bar{\nu}}(\Gamma^+,\Gamma^-)$ of branched covers of the orbit cylinder with one special point and $(s_1)$ in the tautological line bundles $\LL_1$ over all moduli spaces $\CM_{1,1}= \CM_1(\Gamma)$ of branched covers of the complex plane with one special point, we as above can choose coherent collection of sections $(\tilde{s})$ connecting $(s)$ and $(s_1)$ by requiring that over every moduli space $\widetilde{\IM}_1^{\bar{\nu}}$ the section $\tilde{s}$ agrees with the section $s$ over $\CM_1^{\bar{\nu}}$, with the pullback $\pi_1^*s_1$ over  
$\CM_{1,1}\times\CM_{2,0}^{\bar{\nu}}$ and with the pullbacks $\pi_1^*\tilde{s}_1$ over $\widetilde{\IM}_{1,1}^{\bar{\nu}_1}\times\CM_{2,0}^{\bar{\nu}_2}$ and $\pi_2^*\tilde{s}_2$ over $\CM_{1,0}^{\bar{\nu}_1}\times \widetilde{\IM}_{2,1}^{\bar{\nu}_2}$, respectively. \\

Proceeding by induction it then follows that the regular descendant moduli space $\widetilde{\IM}_1^{\bar{\nu},j}$ has codimension-one boundary components $\CM_1^{\bar{\nu},j}$ and $\CM^j_{1,1}\times\CM_{2,0}^{\bar{\nu}}$ as well as $\widetilde{\IM}_{1,1}^{\bar{\nu}_1,j}\times\CM_{2,0}^{\bar{\nu}_2}$ and $\CM_{1,0}^{\bar{\nu}_1}\times \widetilde{\IM}_{2,1}^{\bar{\nu}_2,j}$, respectively. Since we have that $\#\CM_{1,0}^{\bar{\nu}_1}=\#\CM_{2,0}^{\bar{\nu}_2}=0$ by the result in [F] it hence follows that 
\begin{equation*} \#\CM_1^{\bar{\nu},j}=\#\CM^j_{1,1}\cdot\#\CM_{2,0}^{\bar{\nu}} \end{equation*}      
which finally proves the decendants-branching correspondence on the level of coherent collections of sections in obstruction bundles and tautological line bundles. \\

Note in particular that using this stretching process we were able to separate the transversality problem from the problem of defining gravitational descendants. Since the moduli space $\CM_{1,1}=\CM_1(\Gamma)$ is independent of the chosen Reeb orbit and agrees with the moduli space obtained from the degeneration process in relative Gromov-Witten theory, it follows precisely like in the circle bundle case described above that the count of elements in the descendant moduli space $\CM^j_{1,1}$ is independent of the chosen coherent collection of sections and agrees with the integral of the $j$.th power of the first Chern class over $\CM_{1,1}$. $\qed$ 
 
\subsection{Branched covers of trivial half-cylinders}
In the case when the contact manifold $V$ is the unit cotangent bundle $S^*Q$ of a Riemannian manifold $Q$, Cieliebak and Latschev have shown in [CL] that, when suitably interpreted, the symplectic field theory of $V=S^*Q$ without differential forms and gravitational descendants agrees with the string topology of $Q$. The required isomorphism is established by studying punctured holomorphic curves in $T^*Q$ with boundary on the Lagrangian $Q\subset T^*Q$. For this they equip $T^*Q$ with an almost complex structure $\Ju$ such that $(T^*Q,\Ju)$ is an almost complex manifold with one positive cylindrical end $(\IR^+\times S^*Q,\Ju)$. \\

After showing that the contact area of holomorphic curve is given as differences of the sums of the actions of the Reeb orbits in $S^*Q$ and the sum of the lenghts of the boundary components on $Q$, they use the natural filtration by action on symplectic field theory and by length on string topology to show that the morphism has the form of a unitriangular matrix. The entries on the diagonal count cylinders with zero contact area, which are precisely the trivial half-cylinders in $T^*Q$ connecting the geodesic $\bar{\gamma}$ on $Q$ with the corresponding Reeb orbit $\gamma$ in $S^*Q$.  On the other hand, since orbit curves are characterized by the fact that they have zero contact area, it hence directly follows from their proof that there exists a version of their isomorphism statement for the symplectic field theory of a closed Reeb orbit $\gamma$ by studying branched covers over the trivial half-cylinder connecting $\bar{\gamma}$ and $\gamma$. \\

For this let us first recall some definitions from [CL]. Let $\AA^0$ be the graded commutative subalgebra of $\WW$ of polynomials in the variables $q_{\gamma}$, where, following our notation from before, the subscript 0 indicates that no $t$-variables are involved. The Hamiltonian $\IH^0\in\hbar^{-1}\WW^0$ defines a differential operator $\ID^0_{\SFT}:=\overrightarrow{\IH^0}:\AA^0[[\hbar]]\to\AA^0[[\hbar]]$ via the replacements 
\begin{equation*} p_{\gamma}\mapsto\kappa_{\gamma}\hbar\overrightarrow{\frac{\del}{\del q_{\gamma}}}. \end{equation*} 
The resulting pair $(\AA^0[[\hbar]],\ID^0_{\SFT})$ has then the structure of a $\BV_{\infty}$-algebra, in particular, $\ID^0_{\SFT}\circ\ID^0_{\SFT}=0$. Reversely, given a $\BV_{\infty}$-algebra $(\AA^0[[\hbar]],\ID^0)$ where $\AA^0$ is a space of polynomials in variables $q$, it follows, see [CL], that $\ID^0: \AA^0[[\hbar]]\to\AA^0[[\hbar]]$ is a true differential operator. In particular, we naturally get a Weyl algebra $\WW^0$ with distinguished element $\IH^0\in\hbar^{-1}\WW$ satisfying $[\IH^0,\IH^0]=0$ by introducing for each $q$-variable a dualizing $p$-variable, considering the natural commutator relation and using the replacement for $p$-variables from above. \\
 
As already mentioned above, in [CL] it is shown that the $\BV_{\infty}$-algebra $(\AA^0[[\hbar]],\ID_{\SFT})$ representing the symplectic field theory of $S^*Q$ is isomorphic to a $\BV_{\infty}$-algebra $(\CC^0[[\hbar]],\ID^0_{\str})$ constructed from the string topology of $Q$, where $\CC^0$ is a space of chains in the string space $\Sigma=\Sigma Q$ of $Q$. The differential is given by $\ID^0_{\str}=\del + \Delta + \hbar\nabla: \CC^0[[\hbar]]\to\CC^0[[\hbar]]$, where $\del$ is the singular boundary operator and $\nabla$, $\Delta$ are defined using the string bracket and cobracket operations of Chas and Sullivan. The $\BV_{\infty}$-isomorphism $\overrightarrow{\IL^0}$ is defined using the potential of $(T^*Q,Q)$,
\begin{equation*} \IL^0 = \sum_{g,s^-,\Gamma} \ev_{g,s^-}(\Gamma) p^{\Gamma}\hbar^{g-1} \end{equation*} 
using the evaluation cycles $\ev_{g,s^-}(\Gamma)=\ev=(\ev^1,...,\ev^{s^-}): \CM_{g,s^-}(\Gamma)\to\Sigma Q\times ...\times \Sigma Q$ ($s^-$-times) starting from the moduli space of holomorphic curves in $T^*Q$ with positive asymptotics $\Gamma$, genus $g$ and $s^-$ boundary components on $Q$. \\

Now for moving from the symplectic field theory of $S^*Q$ to the symplectic field theory of a closed Reeb orbit $\gamma$ in $S^*Q$, we obviously just have to replace $(\AA^0[[\hbar]],\ID^0_{\SFT}=\overrightarrow{\IH^0})$ by the $\BV_{\infty}$-algebra $(\AA^0_{\gamma}[[\hbar]],\ID^{0,\gamma}_{\SFT}=\overrightarrow{\IH^0_{\gamma}})$ generated only by the $q$-variables representing the multiples of the fixed orbit $\gamma$. Furthermore the potential $\IL^0$ of $(T^*Q,Q)$ is now replaced by the potential $\IL^0_{\gamma,\bar{\gamma}}$ counting branched covers of the trivial half-cylinder connecting $\bar{\gamma}$ in $Q$ and $\gamma$ in $S^*Q$, which defines a $\BV_{\infty}$-isomorphism from $(\AA^0_{\gamma}[[\hbar]],\ID^{0,\gamma}_{\SFT})$ to a $\BV_{\infty}$-algebra $(\CC^0_{\gamma}[[\hbar]],\ID^{0,\bar{\gamma}}_{\str})$. \\

Assigning as for the Reeb orbits formal $q$-variables to multiples of the underlying closed geodesic $\bar{\gamma}$, the potential $\IL^0_{\gamma,\bar{\gamma}}$ is defined by 
\begin{equation*} 
 \IL^0_{\gamma,\bar{\gamma}} = \sum_{g,\bar{\Gamma},\Gamma} \#\CM_g(\Gamma,\bar{\Gamma}) p^{\Gamma} q^{\bar{\Gamma}}\hbar^{g-1} 
\end{equation*} 
summing over all moduli spaces $\CM_g(\Gamma,\bar{\Gamma})$ of branched covers of the trivial half-cylinder with Fredholm index zero. Note that it follows from the area estimate from above for curves in $T^*Q$ with boundary on $Q$ in terms of action of the Reeb orbits and length of the boundary component that, assuming enough transversality, the moduli space $\CM_g(\Gamma,\bar{\Gamma})$ agrees with the preimage of the product stable manifold 
\begin{equation*} W^+(\bar{\Gamma})=W^+(\bar{\gamma}^{n_1})\times ...\times W^+(\bar{\gamma}^{n_{s^-}}) \subset \Sigma Q\times ...\times\Sigma Q \end{equation*} 
of the energy functional $\EE: \Sigma Q\to\IR$ on the string space under the evaluation $\ev_{g,s^-}: \CM_{g,s^-}(\Gamma)\to\Sigma Q\times ...\times \Sigma Q$
\begin{equation*} \CM_g(\Gamma,\bar{\Gamma}) = \ev_{g,s^-}(\Gamma)^{-1}(W^+(\bar{\Gamma}))\subset \CM_{g,s^-}(\Gamma). \end{equation*}

Now the $\BV_{\infty}$-algebra $(\CC^0[[\hbar]],\ID^0_{\str})$ is replaced by the $\BV_{\infty}$-algebra $(\CC^0_{\bar{\gamma}}[[\hbar]],\ID^{0,\bar{\gamma}}_{\str})$ of polynomials in the $q$-variables assigned to multiples of $\bar{\gamma}$. Since this algebra is now indeed an algebra of polynomials, we have seen above that we assign to $(\CC^0_{\bar{\gamma}}[[\hbar]],\ID^{0,\bar{\gamma}}_{\str})$ again a Weyl algebra 
$\WW^0_{\bar{\gamma}}$ with bracket $[\cdot,\cdot]$ generated by $p$- and $q$-variables assigned to multiples of $\bar{\gamma}$ together with a distinguished element $\IG^0_{\bar{\gamma}}\in\hbar^{-1}\WW^0_{\bar{\gamma}}$ satisfying $[\IG^0_{\bar{\gamma}},\IG^0_{\bar{\gamma}}] = 0$. Since $\BV_{\infty}$-algebras $(\AA^0_{\gamma}[[\hbar]],\ID^{0,\gamma}_{\SFT}=\overrightarrow{\IH^0_{\gamma}})$, $(\CC^0_{\bar{\gamma}}[[\hbar]],\ID^{0,\bar{\gamma}}_{\str}=\overrightarrow{\IG^0_{\bar{\gamma}}})$ determine the Weyl algebras with Hamiltonians $(\WW^0_{\gamma},\IH^0_{\gamma})$, $(\WW^0_{\bar{\gamma}},\IG^0_{\bar{\gamma}})$ and vice versa, it follows that the $\BV_{\infty}$-isomorphism given by the potential $\IL^0_{\gamma,\bar{\gamma}}$ indeed leads to an isomorphism of the structures defined by $(\WW^0_{\gamma},\IH^0_{\gamma})$ and  $(\WW^0_{\bar{\gamma}},\IG^0_{\bar{\gamma}})$: \\ 

Indeed, let $\DD^0_{\gamma,\bar{\gamma}}$ be the space of formal power series in the $p$-variables for multiples of $\gamma$ and  
$\hbar$ with coefficients which are polynomials in the $q$-variables assigned to multiples of $\bar{\gamma}$. Then it follows that 
$\IL^0_{\gamma,\bar{\gamma}}$ is an element of $\hbar^{-1}\DD^0_{\gamma,\bar{\gamma}}$ satisfying the master equation
\begin{equation*} 
  e^{\IL^0_{\gamma,\bar{\gamma}}}\overleftarrow{\IH^0_{\gamma}} - \overrightarrow{\IG^0_{\bar{\gamma}}}e^{\IL^0_{\gamma,\bar{\gamma}}} = 0. 
\end{equation*}
In particular, it follows in the notation of section 1 that the map
\begin{equation*} 
 (\IL^{0,+}_{\gamma,\bar{\gamma}})_*^{-1}\circ (\IL^{0,-}_{\gamma,\bar{\gamma}})_*: 
  H_*(\hbar^{-1}\WW^0_{\bar{\gamma}},D^{0,\bar{\gamma}}_{\str})\to H_*(\hbar^{-1}\WW^0_{\gamma},D^{0,\gamma}_{\SFT}) 
\end{equation*}
is an isomorphism of Weyl algebras. \\

In order to understand $D^{0,\bar{\gamma}}_{\str}$, recall that the differential in the string topology was given by $\ID^0_{\str}=\del + \Delta + \hbar\nabla: \CC^0[[\hbar]]\to\CC^0[[\hbar]]$, where $\nabla$ is defined using the string bracket and $\Delta$ using the string cobracket operations defined by Chas and Sullivan. While the singular boundary $\del$ does not appear as we restrict ourselves to zero-dimensional moduli spaces, we expect to get contributions of the string bracket and string cobracket to $\IG^0_{\bar{\gamma}}$, where we claim that the string bracket restricts to the operation of concatenating two multiples $\bar{\gamma}^{n_1}$, $\bar{\gamma}^{n_2}$ to the multiple $\bar{\gamma}^{n_1+n_2}$ of $\bar{\gamma}$, while the string cobracket corresponds to splitting up the multiple $\bar{\gamma}^{n_1+n_2}$ again into $\bar{\gamma}^{n_1}$, $\bar{\gamma}^{n_2}$. \\

In order to see this note that the compactification of the moduli spaces of branched covers of the trivial half-cylinder counted in the potential $\IL^0_{\gamma,\bar{\gamma}}$ can be entirely understood in terms branch points of the branched covering map. While branch points moving the infinite end lead to appearance of $\IH^0_{\gamma}$ in the master equation, the Hamiltonian $\IG^0_{\bar{\gamma}}$ describes what happens if branch points are moving through the boundary of the branched cover, which itself sits over the boundary of the half-cylinder. The important observation is now that for the codimension-one boundary of the moduli space we only have to consider the case where a {\it single} branch point is leaving the branched cover through the boundary.  In order to see that this is described by the concatenation and splitting operations of the multiples of $\bar{\gamma}$, observe that the case when a branch point sits in the boundary of the branched cover is equivalent to the fact that the boundary of $\IR^+\times S^1$ is a critical level set of the branching map followed by the projection to the first factor. Observe that the branch point may leave the branched cover through any point of its boundary, which itself is diffeomorphic to (a number of copies of) the circle. Note that this corresponds to the fact that the concatenation and splitting operation may take place anywhere over any point on $\bar{\gamma}$. It follows that we always get an one-dimensional family of configurations. \\
  
Before we continue, we want to restrict ourselves as before to the rational case. In particular, there exists a version of the above isomorphism, given by counting rational branched covers of the trivial half-cylinder, which relates the rational symplectic field theory $H_*(\PP^0_{\gamma},d^0_{\gamma})$ of $\gamma$ with $H_*(\PP^0_{\bar{\gamma}},d^0_{\bar{\gamma}})$, where $d^0_{\bar{\gamma}}=\{\Ig^0_{\bar{\gamma}},\cdot\}: \PP^0_{\bar{\gamma}}\to\PP^0_{\bar{\gamma}}$ and $\IG^0_{\bar{\gamma}}=\hbar^{-1}(\Ig^0_{\bar{\gamma}}+o(\hbar))$. Before we discuss the rational Hamiltonian $\Ig^0_{\bar{\gamma}}\in\PP^0_{\bar{\gamma}}$, recall that it was shown in [F2] that $\Ih^0_{\gamma}=0$. Note that we have we indeed have not considered additional marked points so far. In particular, it follows from the above isomorphism that also $\Ig^0_{\bar{\gamma}}$ has to vanish. \\

Since we have seen above that for $\IG^0_{\bar{\gamma}}$ and hence for $\Ig^0_{\bar{\gamma}}$ we always get one-dimensional sets of configurations, the vanishing of 
$\Ig^0_{\bar{\gamma}}$ seems to follow from a stupid dimension argument. On the other hand, recall that we have shown in [F2] that the corresponding statement for $\Ih^0_{\gamma}$ does {\it not} simply follow from a symmetry argument but indeed requires a careful study of sections in obstruction bundles in order to find compact perturbations making the Cauchy-Riemann operator transversal to the zero section. With the work in [F2] it is clear that the same transversality problem should continue to hold for branched covers of trivial half-cylinders. In the next section it will turn out that, like on the symplectic field theory side, also on the string side we are working in a highly degenerate situation, so that the transversality requirement is usually not fulfilled. \\

\subsection{Obstruction bundles and transversality}
In order to solve the transversality problem we follow the author's paper [F2] in employing finite-dimensional obstruction bundles over the nonregular configuration spaces. Here is a sketch of the main points. \\

For this let $\Si$ denote a (possibly disconnected) punctured Riemann surface with boundary of genus zero with $s^+$ punctures $z_1,...,z_{s^+}$ and $s^-$ boundary circles $C_1,...,C_{s^-}$ and fix two ordered sets $\Gamma=(\gamma^{n_1^+},...,\gamma^{n_{s^+}^+})$, $\bar{\Gamma}=(\bar{\gamma}^{n_1^-},...,\bar{\gamma}^{n_{s^-}^-})$ of iterates of $\gamma$, $\bar{\gamma}$, respectively. Let $(\xi=TT^*Q/T(\IR^+_0\times S^1),\Ju_{\xi})$ denote the complex normal bundle to the trivial half-cylinder $(\IR^+_0\times S^1,\{0\}\times S^1)\hookrightarrow (T^*Q,Q)$ as defined in [CL], which over the boundary $\{0\}\times S^1\cong\bar{\gamma}\subset Q$ has the property that $\xi\cap TQ$ agrees with the normal bundle $N$ to the geodesic $\bar{\gamma}$ in $Q$. Note that the tangent space $TW^+(\bar{\gamma}^n)$ to the stable manifold of the energy functional in the critical point $\bar{\gamma}^n$ can identified with a subspace of the space of normal deformations $C^0((\bar{\gamma}^n)^*N)$. \\

Given a branched covering $h:(\Si,\del\Si)\to(\IR^+_0\times S^1,\{0\}\times S^1)$ of the trivial half-cylinder, for $p>2$ let $H^{1,p}(h^*\xi)\subset C^0(h^*\xi)$ denote the space of $H^{1,p}$-sections in $h^*\xi$ which over every boundary component $C_k\subset\del\Si$ restrict to a section in $C^0((\bar{\gamma}^{n_k^-})^*N)$. Furthermore we will consider the subspace $H^{1,p}_{\bar{\Gamma}}(h^*\xi)\subset H^{1,p}(h^*\xi)$ consisting of all sections in $h^*\xi$, which over every boundary circle $C_k$ restrict to sections in the subspace $TW^+(\bar{\gamma}^{n_k^-})\subset C^0((\bar{\gamma}^{n_k^-})^*N)$. While the latter Sobolev spaces describe the normal deformations of the branched covering, we introduce similar as in [F2] for sufficiently small $d>0$ a Sobolev space with asymptotic weights $H^{1,p,d}_{\cst}(\Si,\IC)$ in order to keep track of tangential deformations, where, additionally to the definitions in [F2], we impose the natural constraint that the function is real-valued over the boundary. In the same way we define the Banach spaces $L^p(\Lambda^{(0,1)}\Si\otimes_{j,\Ju_{\xi}} h^*\xi)$ and $L^{p,d}(\Lambda^{(0,1)}\Si\otimes_{j,i}\IC)$. Further we denote by $\IM_{0,s^-,s^+}$ the moduli space of Riemann surfaces with $s^-$ boundary circles, $s^+$ punctures and genus zero. \\     

Following [F2], [BM] for the general case and [W] for the case with boundary, there exists a Banach space bundle $\EE$ over a Banach manifold of maps $\BB$ in which the Cauchy-Riemann operator $\CR_J$ extends to a smooth section. In our special case it follows as in [F2] that the fibre is given by  
\begin{equation*} 
 \EE_{h,j} = L^{p,d}(\Lambda^{0,1}\Si\otimes_{j,i}\IC) \oplus L^p(\Lambda^{0,1}\Si\otimes_{j,\Ju_{\xi}}h^*\xi),
\end{equation*}
while the tangent space to the Banach manifold of maps $\BB= \BB_{0,s^-}(\Gamma)$ at $(h,j) \in \IM = \IM_{0,s^-}(\Gamma)$ is given by   
\begin{equation*}
  T_{h,j}\BB = H^{1,p,d}_{\cst}(\Si,\IC)\oplus H^{1,p}(h^*\xi)\oplus T_j\IM_{0,n}.
\end{equation*} 

It follows that the linearization $D_{h,j}$ of the Cauchy-Riemann operator $\CR_{\Ju}$ is a linear map from $T_{h,j}\BB$ to $\EE_{h,j}$, which is surjective in the case when transversality for $\CR_{\Ju}$ is satisfied. In this case it follows from the implicit function theorem that $\ker D_{h,j}=T_{h,j}\IM$. In order to prove that the dimension of the desired moduli space $\IM_{\bar{\Gamma}}=\IM(\Gamma,\bar{\Gamma})=\ev^{-1}(W^+(\bar{\Gamma}))\subset\IM(\Gamma)$ agrees with the virtual dimension expected by the Fredholm index, it remains to prove that the evaluation map $\ev:\IM\to\Sigma Q^{s^-}$ is transversal to the product stable manifold $W^+(\bar{\Gamma})$. \\

In order to deal with this additional transversality problem, we introduce the Banach submanifold of maps $\BB_{\bar{\Gamma}}=\ev^{-1}(W^+(\bar{\Gamma}))\subset\BB$ with tangent space  
\begin{eqnarray*}
  T_{h,j}\BB_{\bar{\Gamma}} &=& H^{1,p,d}_{\cst}(\Si,\IC)\oplus H^{1,p}_{\bar{\Gamma}}(h^*\xi)\oplus T_j\IM_{0,n} \\
                            &=&\{v\in T_{h,j}\BB: v|_{\del\Si}\in TW^+(\bar{\Gamma})\} 
\end{eqnarray*} 
and view the Cauchy-Riemann operator as a smooth section in $\EE\to\BB_{\bar{\Gamma}}$. Then we have the following nice transversality lemma. \\
\\
{\bf Lemma 2.6:} {\it Assume that $D_{h,j}: T_{h,j}\BB_{\bar{\Gamma}} \to \EE_{h,j}$ is surjective. Then the linearization of the evaluation map 
$d_{h,j}\ev: T_{h,j}\IM\to TW^-(\bar{\Gamma})=C^0(\bar{\Gamma}^*N)/TW^+(\bar{\Gamma})$ is surjective.} \\
\\
{\it Proof:} Given $v_0\in TW^-(\bar{\Gamma})$, choose $\tilde{v}\in T_{h,j}\BB$ such that $d_{h,j}\ev\cdot\tilde{v}=v_0$. On the other hand, since 
$D_{h,j}: T_{h,j}\BB_{\bar{\Gamma}} \to \EE_{h,j}$ is onto, we can find $v\in T_{h,j}\BB_{\bar{\Gamma}}$ with $D_{h,j}v=D_{h,j}\tilde{v}$, that is, 
$\tilde{v}-v\in\ker D_{h,j}= T_{h,j}\IM$. On the other hand, since $d_{h,j}\ev\cdot v \in TW^+(\bar{\Gamma})$ for all $v\in T_{h,j}\BB_{\bar{\Gamma}}$ by definition, we have $d_{h,j}\ev\cdot(\tilde{v}-v) = d_{h,j}\ev\cdot\tilde{v} = v_0$ and the claim follows. $\qed$ \\

We have seen that, instead of requiring transversality for the Cauchy-Riemann operator in the Banach space bundle over $\BB$ and geometric transversality for the evaluation map, it suffices to require transversality for the Cauchy-Riemann operator in the Banach space bundle over the smaller Banach manifold $\BB_{\bar{\Gamma}}$. Along the same lines as for proposition 2.1 in [F2] it can be shown that the linearized Cauchy-Riemann operator is of the form 
\begin{eqnarray*}
 &D_{h,j}:& H^{1,p,d}_{\cst}(\Si,\IC)\oplus H^{1,p}_{\bar{\Gamma}}(h^*\xi)\oplus T_j\IM_{0,n} \\
 &&\to L^{p,d}(\Lambda^{0,1}\Si\otimes_{j,i}\IC) \oplus L^p(\Lambda^{0,1}\Si\otimes_{j,\Ju_{\xi}}h^*\xi),\\
 && D_{h,j} \cdot (v_1,v_2,y) = (\CR v_1+ D_j y, D_h^{\xi} v_2),
\end{eqnarray*}
where $\CR: H^{1,p,d}_{\cst}(\Si,\IC) \to L^{p,d}(\Lambda^{0,1}\Si\otimes_{j,i}\IC)$ is the standard 
Cauchy-Riemann operator, $D_h^{\xi}: H^{1,p}(h^*\xi) \to L^p(\Lambda^{0,1}\Si\otimes_{j,\Ju_{\xi}}h^*\xi)$ describes the linearization of $\CR_J$ in the direction of 
$\xi \subset TT^*Q$ and $D_j: T_j\IM_{0,n} \to L^{p,d}(T^*\Si\otimes_{j,i}\IC)$ describes the variation of $\CR_J$ with $j\in\IM_{0,n}$. \\

In [F2] we have shown that for branched covers of orbit cylinders the cokernels of the linearizations of the Cauchy-Riemann operator have the same dimension for every branched cover and hence fit together to give a smooth vector bundle over the nonregular moduli space of branched covers, so that we can prove transversality without waiting for the completion of the polyfold project of Hofer, Wysocki and Zehnder. The following proposition, proved in complete analogy, outlines that this still holds true for branched covers of trivial half-cylinders.   \\
\\
{\bf Proposition 2.7:} {\it The cokernels of the linearizations of the Cauchy-Riemann operator fit together to give a smooth finite-dimensional vector bundle over the moduli space of branched covers of the half-cylinder.} \\
\\
{\it Proof:} As in [F2] this result relies on the transversality of the standard Cauchy-Riemann operator and the super-rigidity of the trivial half-cylinder 
\begin{equation*} \coker\CR=\{0\}\;\; \textrm{and}\;\; \ker D_h^{\xi} =\{0\}, \end{equation*}
where the second statement is now just a linearized version of lemma 7.2 in [CL] which states that, as for orbit cylinders in the symplectizations, the branched covers of the trivial half-cylinder are characterized by the fact that they carry no energy in the sense that the action of Reeb orbits above agrees with the lenghts of the closed geodesics below. $\qed$ \\

It remains to study the extension $\overline{\Coker}\CR_{\Ju}$ of the cokernel bundle $\Coker\CR_{\Ju}$ to the compactified moduli space. For this recall that the components of the codimension-one-boundary of the nonregular moduli space $\CM=\CM_{\bar{\Gamma}}$ of branched covers of the half-cylinder are either of the form $\CM_1\times\CM_2$, where $\CM_1=\overline{\IM_1(\Gamma_1^+,\Gamma_1^-)/\IR}$, $\CM_2=\overline{\IM_2(\Gamma_2,\bar{\Gamma}_2)}$ are nonregular compactified moduli spaces of branched covers of the orbit cylinder or of the trivial half-cylinder, respectively, or of the form $\CM_0\times S^1$, where $\CM_0=\overline{\IM_0(\Gamma,\bar{\Gamma}_0)}$ is again a nonregular compactified moduli space of branched covers of the trivial half-cylinder while $S^1$ refers to the concatenation or splitting locus, which agrees with the locus where the single branch point is leaving the branched covering through the boundary. Note that for $\bar{\Gamma}=(\bar{\gamma}^{n_1},...,\bar{\gamma}^{n_{s^-}})$ the ordered set $\bar{\Gamma}_0$ is either of the form 
\begin{eqnarray*} 
\bar{\Gamma}_0&=&(\bar{\gamma}^{n_1},...,\bar{\gamma}^{n_{k-1}},\bar{\gamma}^{n_k^1},\bar{\gamma}^{n_k^2},\bar{\gamma}^{n_{k+1}},...,\bar{\gamma}^{n_{s^-}}) 
\;\;\textrm{or}\\ \bar{\Gamma}_0&=&(\bar{\gamma}^{n_1},...,\bar{\gamma}^{n_{k-1}},\bar{\gamma}^{n_k+n_{k+1}},\bar{\gamma}^{n_{k+2}},...,\bar{\gamma}^{n_{s^-}}),
\end{eqnarray*}
corresponding to concatenating $\bar{\gamma}^{n_k^1}$ and $\bar{\gamma}^{n_k^1}$ to get $\bar{\gamma}^{n_k}$ ($n_k^1+n_k^2=n_k$) or the splitting of $\bar{\gamma}^{n_k+n_{k+1}}$ to get $\bar{\gamma}^{n_k}$ and $\bar{\gamma}^{n_{k+1}}$. Restricting to the concatenation case, recall that the chosen special point on the simple closed Reeb orbit determines a special point on the underlying simple geodesic and that we may assume that every holomorphic curve comes equipped with asymptotic markers in the sense of [EGH] not only on the cylindrical ends but also on the boundary circles. In particular, for the concatenation and splitting processes we may assume that all multiply-covered geodesics come equipped with a parametrization by $S^1$. Denoting by $t_1, t_2\in S^1$ the points on $\bar{\gamma}^{n_k^1}$, $\bar{\gamma}^{n_k^1}$, where we want to concatenate the two multiply-covered geodesics to get the multiply-covered geodesic $\bar{\gamma}^{n_k^1+n_k^2}$, we see that the coordinates must satisfy $n_k^1 t_1=n_k^2 t_2$ in order to represent the same point on the underlying simple geodesic, so that the configuration space agrees with $S^1$ by setting $t_1=n_k^2 t$, $t_2=n_k^1 t$ for $t\in S^1$.  \\

While it directly follows from [F2] that over the boundary components $\CM_1\times\CM_2\subset\CM$ the extended cokernel bundle $\overline{\Coker}\CR_{\Ju}$ is of the form 
\begin{equation*} \overline{\Coker}\CR_{\Ju}|_{\CM_1\times\CM_2} = \pi_1^*\overline{\Coker}^1\CR_{\Ju}\oplus \pi_2^*\overline{\Coker}^2\CR_{\Ju}, \end{equation*}
where $\overline{\Coker}^1\CR_{\Ju}$, $\overline{\Coker}^2\CR_{\Ju}$ denote the (extended) cokernel bundles over $\CM_1$, $\CM_2$, respectively, it remains to study the  cokernel bundle over the boundary components $\CM_0\times S^1$. \\
\\
{\bf Proposition 2.8:} {\it Over the boundary components $\CM_0\times S^1\subset \CM$ the extended cokernel bundle $\overline{\Coker}\CR_{\Ju}$ is also of product form,}
\begin{equation*}  \overline{\Coker}\CR_{\Ju}|_{\CM_0\times S^1} \;=\; \pi_1^*\overline{\Coker}^0\CR_{\Ju}\oplus \pi_2^*\Delta, \end{equation*}
{\it where $\overline{\Coker}^0\CR_{\Ju}$ denotes the (extended) cokernel bundle over the moduli space $\CM_0$ and $\Delta$ is a vector bundle over $S^1$ which is determined by the tangent spaces to the stable manifolds of the multiply-covered closed geodesics involved into the concatenation or splitting process.}\\  
\\
{\it Proof:} Still restricting to the concatenation case, let $\Si_0=\Si_{01}\cup\Si_{02}$ denote the disconnected Riemann surface of genus zero with $s^+$ punctures and $s^-+1$ boundary circles $C_1,...,C_k^1,C_k^2,...,C_{s^-}$, where we assume that $\del\Si_{01}=C_1\cup...\cup C_k^1$ and $\del\Si_{02}= C_k^2,...,C_{s^-}$. As before we know that the tangent spaces to the corresponding Banach manifolds of maps $\BB^0$, $\BB^0_{\bar{\Gamma}_0}$ at a branched covering $(h_0,j_0): (\Si_0,\del\Si_0)\to(\IR^+_0\times S^1,\{0\}\times S^1)$ are given by 
\begin{eqnarray*}
  T_{h_0,j_0}\BB^0 &=& H^{1,p,d}_{\cst}(\Si_0,\IC)\oplus H^{1,p}(h_0^*\xi)\oplus T_{j_0}\IM_{0,n},\\
  T_{h_0,j_0}\BB^0_{\bar{\Gamma}_0} &=& H^{1,p,d}_{\cst}(\Si_0,\IC)\oplus H^{1,p}_{\bar{\Gamma}_0}(h_0^*\xi)\oplus T_{j_0}\IM_{0,n} \\
                                  &=&\{v\in T_{h_0,j_0}\BB^0: v|_{\del\Si_0}\in TW^+(\bar{\Gamma}_0)\} 
\end{eqnarray*} 
while the fibre of the corresponding Banach space bundle is given by
\begin{equation*} 
 \EE^0_{h_0,j_0} = L^{p,d}(\Lambda^{0,1}\Si_0\otimes_{j_0,i}\IC) \oplus L^p(\Lambda^{0,1}\Si_0\otimes_{j_0,\Ju_{\xi}}h_0^*\xi).
\end{equation*}
For $(h_0,j_0,t)\in\CM_0\times S^1$ we further introduce the Banach manifold of maps $\BB^*_{\bar{\Gamma}}\subset\BB^*\subset\BB^0$ which should consist of all branched covers of the trivial half-cylinder in $\BB^0$ for which the boundary circles $C_k^1, C_k^2\cong S^1$ are concatenated at $(t_1,t_2)=(n_k^2 t, n_k^1 t)\in C_k^1\times C_k^2$, to give the singular Riemann surface $\Si_*$ with $s^-$ boundary circles $C_1,...,C_k^1\cup_t C_k^2,...,C_{s^-}$ and we have
\begin{eqnarray*} 
  T_{h_0,j_0,t}\BB^*&=&\{v\in T_{h_0,j_0}\BB^0: v_k^1(n_k^2 t)=v_k^2(n_k^1 t)\;\textrm{for}\; v_k^{1,2}:=v|_{C_k^{1,2}}\} \\
  T_{h_0,j_0,t}\BB^*_{\bar{\Gamma}} &=&\{v\in T_{h_0,j_0,t}\BB^*: v|_{\del\Si_*}\in TW^+(\bar{\Gamma}_0)\}. 
\end{eqnarray*}

The proof of the general gluing theorem in [MDSa] suggests that over $(h_0,j_0,t)\in\CM_0\times S^1\subset\CM$ the extended cokernel bundle $\overline{\Coker}\CR_{\Ju}$ has fibre 
\begin{equation*} 
   (\overline{\Coker}\CR_{\Ju})_{h_0,j_0,t} = \coker D_{h_0,j_0,t},\;\; D_{h_0,j_0,t}: T_{h_0,j_0,t}\BB^*_{\bar{\Gamma}}\to\EE^0_{h_0,j_0}.
\end{equation*}
Before we describe the relation to the cokernel bundle $\overline{\Coker}^0\CR_{\Ju}$ over the first factor $\CM_0$ with fibre
\begin{equation*} 
   (\overline{\Coker}^0\CR_{\Ju})_{h_0,j_0} = \coker D_{h_0,j_0},\;\; D_{h_0,j_0}: T_{h_0,j_0,t}\BB^0_{\bar{\Gamma}_0}\to\EE^0_{h_0,j_0}.
\end{equation*}
observe that we still have 
\begin{eqnarray*} 
\coker D_{h_0,j_0} = \coker D^{\xi}_{h_0},&& D^{\xi}_{h_0}: T^{\xi}_{h_0,j_0}\BB^0_{\bar{\Gamma}_0}\to\EE^{0,\xi}_{h_0,j_0},\\
\coker D_{h_0,j_0,t} = \coker D^{\xi}_{h_0,t},&& D^{\xi}_{h_0,t}: T^{\xi}_{h_0,j_0,t}\BB^*_{\bar{\Gamma}}\to\EE^{0,\xi}_{h_0,j_0},
\end{eqnarray*}
and $\ker D^{\xi}_{h_0}=\ker D^{\xi}_{h_0,t}=\{0\}$, where $T^{\xi}_{h_0,j_0}\BB^0_{\bar{\Gamma}_0}\subset T_{h_0,j_0}\BB^0_{\bar{\Gamma}_0}$, $T^{\xi}_{h_0,j_0,t}\BB^*_{\bar{\Gamma}}\subset T_{h_0,j_0,t}\BB^*_{\bar{\Gamma}}$ and $\EE^{0,\xi}_{h_0,j_0}\subset \EE^0_{h_0,j_0}$ are the subspaces corresponding to normal deformations. \\

Now observing that 
\begin{equation*} 
  TW^+(\bar{\gamma}^{n_k^1})\oplus TW^+(\bar{\gamma}^{n_k^2}) \subset \{v\in TW^+(\bar{\gamma}^{n_k}): v_k^1(n_k^2 t) = v_k^2(n_k^1 t)\}
\end{equation*}
we get from 
\begin{eqnarray*} 
   T^{\xi}_{h_0,j_0}\BB^0_{\bar{\Gamma}_0} =\{v\in T^{\xi}_{h_0,j_0}\BB^0: v|_{\del\Si_0}\in TW^+(\bar{\Gamma}_0)\},\\
   T^{\xi}_{h_0,j_0,t}\BB^0_{\bar{\Gamma}} =\{v\in T^{\xi}_{h_0,j_0,t}\BB^*: v|_{\del\Si_0}\in TW^+(\bar{\Gamma})\}\\
\end{eqnarray*}
that $T^{\xi}_{h_0,j_0}\BB^0_{\bar{\Gamma}_0}\subset T^{\xi}_{h_0,j_0,t}\BB^0_{\bar{\Gamma}}$ with quotient space
\begin{equation*} 
 \frac{T^{\xi}_{h_0,j_0}\BB^0_{\bar{\Gamma}_0}}{T^{\xi}_{h_0,j_0,t}\BB^0_{\bar{\Gamma}}} \;=\; 
 \frac{TW^+(\bar{\gamma}^{n_k^1})\oplus TW^+(\bar{\gamma}^{n_k^2})}{\{v\in TW^+(\bar{\gamma}^{n_k}): v_k^1(n_k^2 t) = v_k^2(n_k^1 t)\}}.
\end{equation*}
On the other hand, since $\ker D^{\xi}_{h_0}=\ker D^{\xi}_{h_0,t}=\{0\}$ we also find for the quotient space that
\begin{equation*} 
 \frac{T^{\xi}_{h_0,j_0}\BB^0_{\bar{\Gamma}_0}}{T^{\xi}_{h_0,j_0,t}\BB^0_{\bar{\Gamma}}} \;=\; 
 \frac{\im D^{\xi}_{h_0}}{\im D^{\xi}_{h_0,t}} \;=\; \frac{\coker D^{\xi}_{h_0,t}}{\coker D^{\xi}_{h_0}},
\end{equation*}
where the last equality follows from the fact that $D^{\xi}_{h_0}$ and $D^{\xi}_{h_0,t}$ both map to the same Banach space $\EE^{0,\xi}_{h_0}$. \\

Defining an obstruction bundle $\Delta$ over $S^1$ by setting 
\begin{equation*}
 \Delta_t = \frac{TW^+(\bar{\gamma}^{n_k^1})\oplus TW^+(\bar{\gamma}^{n_k^2})}{\{v\in TW^+(\bar{\gamma}^{n_k}): v_k^1(n_k^2 t) = v_k^2(n_k^1 t)\}}
\end{equation*}
and putting everything together we hence found that 
\begin{equation*} 
  (\overline{\Coker}\CR_{\Ju})_{h_0,j_0,t} \cong (\overline{\Coker}^0\CR_{\Ju})_{h_0,j_0} \oplus \Delta_t,
\end{equation*}
as desired. $\qed$ \\
 
With this we can prove the desired statement about $\Ig^0_{\bar{\gamma}}$. \\
\\
{\bf Corollary 2.9:} {\it We have $\Ig^0_{\bar{\gamma}}=0$.} \\
\\
{\it Proof:} It follows that the obstruction bundle over the one-dimensional configuration space has rank
\begin{equation*} \rank \Delta = \Morse(\bar{\gamma}^{n_k})-\Morse(\bar{\gamma}^{n_k^1})-\Morse(\bar{\gamma}^{n_k^2})+\dim Q-1 \geq 0,
\end{equation*}
where the latter inequality can be verified as in [F2] using the multiple cover index formulas in [Lo]. When by index reasons the configuration is expected to be discrete we get a rank-one obstruction bundle over the boundary of the branched cover, which by orientability reasons must indeed be trivial. $\qed$ \\

On the other hand, we want to emphasize that the proof of $\Ig^0_{\bar{\gamma}}=0$ is much simpler than the proof of $\Ih^0_{\gamma}=0$ in [F2], which has to involve obstruction bundles of arbitrary large rank and uses induction. Besides that our proof in [F2] also holds for Reeb orbits in general contact manifolds, this does not come as surprise. Going back to the symplectic field theory of unit cotangent bundles $S^*Q$, it is already mentioned in [CL] that the SFT differential $\ID^0_{\SFT}=\overrightarrow{\IH^0}:\AA^0[[\hbar]]\to\AA^0[[\hbar]]$ involving all moduli spaces of holomorphic curves in $\IR\times S^*Q$ is much larger than the string differential $\ID^0_{\str}=\del+\Delta+\hbar\nabla:\CC^0[[\hbar]]\to\CC^0[[\hbar]]$, which just involves the singular boundary operator and the string bracket and cobracket operations. \\

\subsection{Additional marked points and gravitational descendants}
We now want to understand the system of commuting operators defined for Reeb orbits by studying moduli spaces of branched covers over the cylinder over $\gamma$ in terms of operations defined for the underlying closed geodesic $\bar{\gamma}$. To this end we have to extend the picture of [CL] used for computing the symplectic field theory of Reeb orbits to include additional marked points on the moduli spaces, integration of differential forms and gravitational descendants. \\

Reintroducing the sequence of formal variables $t_j$, $j\in\IN$, we now consider the graded Weyl algebras $\WW_{\gamma}$, $\WW_{\bar{\gamma}}$ of power series in $\hbar$, the $p$-variables corresponding to multiples of $\gamma$, $\bar{\gamma}$ and $t$-variables with coefficients which are polynomials in the $q$-variables corresponding to multiples of $\gamma$, $\bar{\gamma}$. In the same way we can introduce the graded commutative algebras $\AA_{\gamma}$, $\CC_{\bar{\gamma}}$ of power series in $\hbar$, the $t$-variables with coefficients which are polynomials in the $q$-variables corresponding to multiples of $\gamma$, $\bar{\gamma}$. For the expansion $\IH_{\gamma}=\IH^0_{\gamma}+\sum_j t_j \IH^1_{\gamma,j}+o(t^2)$ of the Hamiltonian from before, we are hence looking for an extended potential $\IL_{\gamma,\bar{\gamma}}$ as well as extended string Hamiltonian $\IG_{\bar{\gamma}}$,  
\begin{eqnarray*}
 \IL_{\gamma,\bar{\gamma}}&=&\IL^0_{\gamma,\bar{\gamma}}+\sum_j t_j \IL^1_{\gamma,\bar{\gamma},j}+o(t^2), \\
 \IG_{\bar{\gamma}}&=&\IG^0_{\bar{\gamma}}+\sum_j t_j \IG^1_{\bar{\gamma},j}+o(t^2),
\end{eqnarray*}
such that $\overrightarrow{\IL_{\gamma,\bar{\gamma}}}: (\AA_{\gamma}[[\hbar]],\overrightarrow{\IH_{\gamma}}) \to (\CC_{\bar{\gamma}}[[\hbar]],\overrightarrow{\IG_{\bar{\gamma}}})$ is an isomorphism of $\BV_{\infty}$-algebras. 
For this we have to prove the extended master equation
\begin{equation*} 
  e^{\IL_{\gamma,\bar{\gamma}}}\overleftarrow{\IH_{\gamma}} - \overrightarrow{\IG_{\bar{\gamma}}}e^{\IL_{\gamma,\bar{\gamma}}} = 0, 
\end{equation*}
while the isomorphism property again follows using the natural filtration given by the $t$-variables. \\

Since we are only interested in the system of commuting operators $\IH^1_{\gamma,j}$, $j\in\IN$, which is defined by counting branched covers of orbit cylinders with at most one additional marked point, we again will only discuss the required compactness statements in the case of one additional marked point. Furthermore we will still just restrict to the rational case. In other words we will prove the following proposition, which is just a reformulation of our theorem from above. \\
\\
{\bf Proposition 2.10:} {\it The system of Poisson-commuting functions $\Ih^1_{\gamma,j}$, $j\in\IN$ on $\PP^0_{\gamma}$ is isomorphic to a system of Poisson-commuting functions $\Ig^1_{\bar{\gamma},j}$, $j\in\IN$ on $\PP^0_{\bar{\gamma}}=\PP^0_{\gamma}$, where for every $j\in\IN$ the descendant Hamiltonian $\Ig^1_{\bar{\gamma},j}$ given by} 
\begin{equation*} 
 \Ig^1_{\bar{\gamma},j} \;=\; \sum \epsilon(\vec{n})\frac{q_{n_1}\cdot ... \cdot q_{n_{j+2}}}{(j+2)!} 
\end{equation*}
{\it where the sum runs over all ordered monomials $q_{n_1}\cdot ... \cdot q_{n_{j+2}}$ with $n_1+...+n_{j+2} = 0$ \textbf{and which are of degree $2(m+j-3)$}. Further $\epsilon(\vec{n})\in\{-1,0,+1\}$ is fixed by a choice of coherent orientations in symplectic field theory and is zero if and only if one of the orbits $\gamma^{n_1},...,\gamma^{n_{j+2}}$ is bad.} \\

{\it Proof:} While the proof seems to require the definition of gravitational descendants for moduli spaces of holomorphic curves not only with punctures but also with boundary, instead of defining them recall that we have shown in the previous subsection 2.2 that the gravitational descendants can be replaced by imposing branching conditions over the special marked point on the orbit cylinder. More precisely, recall the lemma in subsection 2.2 states that we can indeed write each of the Hamiltonians $\Ih^1_{\gamma,j}$ as a weighted sum, 
\begin{equation*} 
 \Ih^1_{\gamma,j} \;=\; \frac{1}{j!}\;\cdot\;\Ih^1_{\gamma,(j)} \;+\;\sum_{|\mu|<j} \rho^0_{j,\mu}\;\cdot\; \Ih^1_{\gamma,\mu}, 
\end{equation*} 
where $\Ih^1_{\gamma,\mu}\in\PP^0_{\gamma}$ counts rational branched covers of the orbit cylinder with $\ell(\mu)$ connected components carrying precisely one additional marked point $z_1,...,z_{\ell(\mu)}$, which are mapped to the special point on the orbit cylinder and $z_i$ is a branch point of order $\mu_i-1$ for all $i=1,...,\ell(\mu)$. \\

While for the invariance statement for gravitational descendants we were studying the compactification of the moduli spaces of holomorphic curves with one additional marked points, it follows from the definition of $\Ih^1_{\gamma,\mu}$ that now it is natural to study the moduli spaces of branched covers of the trivial half-cylinder with $\ell(\mu)$ connected components carrying precisely one additional marked point $z_1,...,z_{\ell(\mu)}$, which are mapped to the special point on the trivial half-cylinder and $z_i$ is a branch point of order $\mu_i-1$ for all $i=1,...,\ell(\mu)$. While for the orbit cylinder the natural $\IR$-action is used to fix not only the $S^1$-coordinate but also the $\IR$-coordinate of the special point, note that, in order to find the branched covers of the orbit cylinder counted in $\Ih^1_{\gamma,\mu}$ in the boundary, for the trivial half-cylinder we still fix $S^1$-coordinate but allow the $\IR$-coordinate to vary in $\IR^+=(0,\infty)$. \\

It follows that besides the boundary phenomena of the moduli spaces of branched covers of the trivial half-cylinder already described above, which can be described as seen above as the moving of branch points to infinity or leaving the branched cover through the boundary, the new boundary phenomena are the moving of the additional marked points to infinity or leaving the branched cover through the boundary, which are equivalent to the moving of the special point to infinity or leaving the half-cylinder through the boundary. In particular, it follows from the latter equivalence that the additional marked points $z_1,...,z_{\ell(\mu)}$ move to infinity or leave the branched cover all at once. While the moving of the additional marked points to infinity, possibly together with other branch points, is counted in $\Ih^1_{\gamma,\mu}$, the corresponding string Hamiltonian $\Ig^1_{\bar{\gamma},\mu}$ should describe what happens if the additional marked points leave the branched cover through the boundary. Provided that we have found $\Ig^1_{\bar{\gamma},\mu}\in\PP^0_{\bar{\gamma}}$ for all branching profiles $\mu$, it then follows from linearity that we obtain the desired Poisson-commuting sequence $\Ig^1_{\bar{\gamma},j}$ by setting 
\begin{equation*} 
 \Ig^1_{\bar{\gamma},j} \;=\; \frac{1}{j!}\;\cdot\;\Ig^1_{\bar{\gamma},(j)} \;+\;\sum_{|\mu|<j} \rho^0_{j,\mu}\;\cdot\; \Ig^1_{\bar{\gamma},\mu}. 
\end{equation*} 

On the other hand, recall that in the computation of $\Ig^0_{\bar{\gamma}}$ we were faced with a transversality problem. While we have shown that the set of configurations counted for $\Ig^0_{\bar{\gamma}}$ is always one-dimensional, one can compute using the Morse indices of the involved multiply-covered geodesics that it happens that the Fredholm index expects the same set to be discrete. In the case when the Fredholm index is right, we have shown that get a obstruction bundle of rank one to cut down the dimension of the configuration space, which is however trivial by orientability. For $\Ig^1_{\bar{\gamma},\mu}$ we now show that the situation is even nicer. \\
\\
{\bf Lemma 2.11:} {\it For every branching condition $\mu$ the set of configurations studied for $\Ig^1_{\bar{\gamma},\mu}$ is already discrete \textbf{before} we add abstract perturbations to the Cauchy-Riemann operator. It follows that, if the Fredholm index is right, there is \textbf{no} obstruction bundle.} \\  

Before we show why this lemma leads to a proof of the above proposition and hence of the theorem, note that when $\bar{\gamma}=Q=S^1$ transversality is always satisfied and hence there are no obstruction bundles at all. On the other hand, note that the above proposition is formulated such that it holds in this case, where we use that $\Ig^1_{S^1,\mu}=\Ih^1_{S^1,\mu}$, which follows from the fact that the (rational) potential $\IL^0_{S^1,S^1}$ ($\Il^0_{S^1,S^1}$) only counts orbit cylinders. \\

In order to see that for an arbitrary closed geodesic $\bar{\gamma}\subset Q$ the lemma proves the proposition and hence the theorem, observe that the Fredholm index is right precisely when it leads to the maximal degree $2(m+j-3)$ from the proposition. Since the configuration space is independent of $\bar{\gamma}$ before perturbing, in this case the lemma tells us that the corresponding configurations counted for $\Ig^1_{\bar{\gamma},\mu}$ indeed agree with the ones counted for $\Ig^1_{S^1,\mu}$, {\it up to sign} determined by a choice of coherent orientations for the moduli spaces as described in [BM]. On the other hand, the results in [BM] show that the bad orbits indeed cancel out. For both statements we refer to the work of Cieliebak and Latschev in order to show that the orientation choices for closed Reeb orbits have a natural translation into orientation choices for to the underlying closed geodesics, that is, their unstable manifolds for the energy functional. In particular, we have, see [CL], that the Reeb orbit $\gamma$ is bad if and only if the unstable manifold of $\bar{\gamma}$ is not orientable. On the other hand, when the Fredholm is not right and hence maximal, we do not get a contribution to $\Ig^1_{\bar{\gamma},\mu}$ by definition. $\qed$\\

Hence it just remains to prove the lemma.\\
\\  
{\it Proof of the lemma:} For simplicity we first prove the statement for $\mu=(2)$. Following the above description of $\Ig^1_{\bar{\gamma,\mu}}$ it follows that $\Ig^1_{\bar{\gamma},(2)}$ describes what happens if the additional marked point, which is a simple branch point, leaves the branched cover through the boundary. While at first this sounds that $\Ig^1_{\bar{\gamma},(2)}$ agrees with $\Ig^0_{\bar{\gamma}}$, note that now the branch point is required to sit over the special point on the boundary of the half-cylinder. Since the $S^1$-coordinate of the special point is fixed, it follows that the branch point can no longer leave the branched cover through every point on the boundary. In particular, while for $\Ig^0_{\bar{\gamma}}$ we obtained a one-dimensional configuration space due to the obvious $S^1$-symmetry, it follows that for the configurations counted in $\Ig^1_{\bar{\gamma},(2)}$ the $S^1$-symmetry is no longer present. Due to the important observation (which we already used to compute $\Ig^0_{\bar{\gamma}}$) that for the codimension-one boundary we can assume that there are no other branch points leaving the boundary at the same time, it follows that the set of configurations is indeed discrete. On the other hand, it is clear that this argument immediately generalizes to all branching profiles $\mu$, since all the $\ell(\mu)$ additional marked points are mapped to the same fixed special point. Together with the observation that the additional marked points $z_1,...,z_{\ell(\mu)}$ leave the branched cover through the boundary all at once when the special point leaves the half-cylinder through the boundary, but again no other branch points by codimension reasons, the corresponding set of configurations stays discrete. $\qed$

\end{document}
